\documentclass[12pt, a4paper]{amsart}

\usepackage[english]{babel}

\usepackage[all]{xy}
\usepackage{amsmath}
\usepackage{amssymb}
\usepackage{amsfonts}
\usepackage{graphicx}
\usepackage[mathscr]{eucal}
\usepackage{enumerate}
\usepackage[latin1]{inputenc}
\usepackage{times}

\theoremstyle{definition}
\newtheorem{definition}{Definition}[section]
\theoremstyle{remark}

\newtheorem{example}[definition]{Example}
\newtheoremstyle{estilo}{\topsep}{\topsep}{\slshape}{}{\bfseries}{.}{ }{}
\theoremstyle{estilo}
\newtheorem{theorem}[definition]{Theorem}
\newtheorem{corollary}[definition]{Corollary}
\newtheorem{lemma}[definition]{Lemma}
\newtheorem{proposition}[definition]{Proposition}
\newtheorem{remark}[definition]{Remark}

\def\mf{\mathfrak}
\def\mb{\mathbb}
\newcommand{\mfa}{\mbox{$\mf {a}$}}
\newcommand{\mfb}{\mbox{$\mf {b}$}}
\def\ot{\otimes}
\newcommand{\gsm}{\mbox{$\blacktriangleright \hspace{-0.7mm}<$}}
\def\cd{\cdot}
\newcommand{\gtl}{\mbox{${\;}$$>\hspace{-0.85mm}\blacktriangleleft$${\;}$}}
\newcommand{\smi}{\mbox{$S^{-1}$}}
\def\r{\rho}

\renewenvironment{proof}{{\noindent\sc Proof\;}}{\\}

\DeclareMathOperator{\Id}{Id}

\def\dirlim{\varinjlim}

\def\O{\Omega}
\def\c{\mathbb C}
\newcommand{\dtext}[1]{\emph{\textbf{#1}}}
\def\vphi{\varphi}
\def\eps{\varepsilon}
\newcommand{\abs}[1]{\left\lvert#1\right\rvert}
\def\M{{\mathcal{M}}}
\def\R{\mathbb R}
\def\n{\mathbb N}
\def\z{\mathbb Z}
\def\l{\lambda}
\def\lto{\longrightarrow}
\def\lmto{\longmapsto}
\newcommand{\imp}[2]{\noindent$\boxed{\textbf{#1}\mathbf{\Rightarrow}\textbf{#2}}\quad$} 
\newcommand{\esc}[1]{\left\langle #1\right\rangle}

\title{On iterated twisted tensor products of algebras}
\thanks{Pascual Jara and Javier López have been
    partially supported by projects MTM2004-08125 and FQM-266 (Junta de Andalucía Research Group).
    Javier López has also been supported by the Spanish MEC FPU-grant AP2003-4340 and the European
    Science Foundation Programme on NONCOMMUTATIVE GEOMETRY (NOG). Florin Panaite and Fred Van Oystaeyen have been
    partially supported by the EC programme
LIEGRITS, RTN 2003, 505078, and by the bilateral project ``New
techniques in Hopf algebras and graded ring theory'', of the Flemish
and Romanian Ministries of Research. Florin Panaite has also been
partially supported by the CEEX programme of the Romanian Ministry
of Education and Research, contract nr. CEx05-D11-11/2005
    }
\author{Pascual Jara Martínez}
\address{Department of Algebra,
University of Granada\\
Avda. Fuentenueva s/n, E-18071, Granada, Spain}
\email{pjara@ugr.es}
\author {Javier López Peña}
\address{Department of Algebra,
University of Granada\\
Avda. Fuentenueva s/n, E-18071, Granada, Spain}
\email{jlopez@ugr.es}
\author{Florin Panaite}
\address{ Institute of Mathematics of the Romanian Academy.\\
PO-Box 1-764, RO-014700 Bucharest (Romania)
}
\email{Florin.Panaite@imar.ro}
\author{ Fred Van Oystaeyen }
\address{Department of Mathematics and Computer Sciences, University of Antwerp.\\ Middelheimlaan 1, B-2020 Antwerp
    (Belgium)}
\email{Francine.Schoeters@ua.ac.be (secretary)}
\date{}

\begin{document}

\newcommand{\Cruce}
    {\xy (0,12)*{}; (8,0)*{} **\crv{(0,4)&(8,6)} \POS?(.5)*{\hole}="x";
         (0,0)*{}; "x" **\crv{(0,4)};
         "x"; (8,12)*{} **\crv{(8,6)};
    \endxy}

\newcommand{\CruceInv}
    {\xy (0,0)*{}; (8,12)*{} **\crv{(0,4)&(8,6)} \POS?(.5)*{\hole}="x";
         (0,12)*{}; "x" **\crv{(0,6)};
         "x"; (8,0)*{} **\crv{(8,4)};
    \endxy}

\newcommand{\CruceDoble}
    {\xy (0,12)*{}; (16,0)*{} **\crv{(0,4)&(16,6)} \POS?(.4)*{\hole}="x";
         (8,12)*{}; "x" **\crv{(8,9)&(0,0)}; "x"; (0,0)*{} **\crv{(0,3)};
         (16,12)*{}; "x" **\crv{(16,9)&(8,6)}; "x"; (8,0)*{} **\crv{(8,3)};
    \endxy}

\newcommand{\CruceDobleInv}
    {\xy (8,12)*{}; (16,0)*{} **\crv{(8,4)&(16,6)} \POS?(.4)*{}="x";
         (0,12)*{}; "x"-(2,0)*{}; **\crv{(0,8)};
         "x"-(2,0)*{}; (10,2)*{}; **\crv{(15,4)};
         (10,2)*{}; (8,0)*{}; **\crv{(8,1)};
         (16,12)*{}; "x"+(2,0); **\crv{(16,8)}; "x"-(2,1.5); (0,0)*{} **\crv{(0,2)};
    \endxy}

\newcommand{\CruceT}[1]
    {\xy (0,12)*{}; (8,0)*{} **\crv{(0,4)&(8,6)} \POS?(.5)*{\hole}="x"+(0,3)*{\scriptstyle{#1}};
         (0,0)*{}; "x" **\crv{(0,4)};
         "x"; (8,12)*{} **\crv{(8,6)};
    \endxy}

\newcommand{\CruceInvT}[1]
    {\xy (0,0)*{}; (8,12)*{} **\crv{(0,4)&(8,6)} \POS?(.5)*{\hole}="x"+(0,3)*{\scriptstyle{#1}};
         (0,12)*{}; "x" **\crv{(0,6)};
         "x"; (8,0)*{} **\crv{(8,4)};
    \endxy}

\newcommand{\Mult}
    {   \xy (0,12)*{}; (8,12)*{} **\crv{(0,3)&(8,3)} \POS?(.5)*{}="x";
            "x"*{}; (4,0)*{} **\crv{(4,4)};
        \endxy
    }

\newcommand{\Unit}[1]
    {   \xy (0,6)*{\scriptscriptstyle{#1}}*\cir<5.85pt>{}; (0,0)*{} **\crv{-};
        \endxy
    }

\newcommand{\CoMult}
    {   \xy (0,0)*{}; (8,0)*{} **\crv{(0,7)&(8,7)} \POS?(.5)*{}="x";
            (4,12)*{}; "x"*{} **\crv{(4,4)};
        \endxy
    }

\newcommand{\CoMultL}[1]
    {   \xy (0,0)*{}; (8,0)*{} **\crv{(0,7)&(8,7)} \POS?(.5)*{}="x";
            (4,12)*{}; "x"*{} **\crv{(4,4)};
            "x"+(4,2)*{\scriptstyle #1};
        \endxy
    }

\newcommand{\MoveRight}
    {   \xy (0,12)*{}; (4,0)*{} **\crv{(0,5)&(4,5)};
        \endxy
    }

\newcommand{\BigMoveRight}
    {   \xy (0,12)*{}; (8,0)*{} **\crv{(0,5)&(8,5)};
        \endxy
    }

\newcommand{\VeryBigMoveRight}
    {   \xy (0,12)*{}; (16,0)*{} **\crv{(0,5)&(16,5)};
        \endxy
    }

\newcommand{\LongBigMoveRight}
    {   \xy (0,24)*{}; (12,0)*{} **\crv{(0,18)&(12,5)};
        \endxy
    }

\newcommand{\LongMoveRight}
    {   \xy (0,24)*{}; (8,0)*{} **\crv{(0,18)&(8,5)};
        \endxy
    }

\newcommand{\MoveLeft}
    {   \xy (4,12)*{}; (0,0)*{} **\crv{(4,5)&(0,5)};
        \endxy
    }

\newcommand{\BigMoveLeft}
    {   \xy (8,12)*{}; (0,0)*{} **\crv{(8,5)&(0,5)};
        \endxy
    }

\newcommand{\LongBigMoveLeft}
    {   \xy (12,24)*{}; (0,0)*{} **\crv{(12,18)&(0,5)};
        \endxy
    }

\newcommand{\Down}{\xy (0,12)*{}; (0,0)*{} **\crv{-}; \endxy}

\newcommand{\DownL}[1]{\xy (0,12)*{}; (0,6)*{\scriptscriptstyle{\vphantom{f}#1}}*+<1.2pt>\frm{-}="x"; **\crv{-};
    "x"; (0,0)*{} **\crv{-}; \endxy}

\newcommand{\DownLb}[1]{\xy (0,12)*{}; (0,6)*{\boxed{\scriptscriptstyle{#1}}}="x"; **\crv{-};
    "x"; (0,0)*{} **\crv{-}; \endxy}


\begin{abstract}
We introduce and study the definition, main properties and
applications of iterated twisted tensor products of algebras,
motivated by the problem of defining a suitable representative for
the product of spaces in noncommutative geometry. We find conditions
for constructing an iterated product of three factors, and prove
that they are enough for building an iterated product of any number
of factors. As an example of the geometrical aspects of our
construction, we show how to construct differential forms and
involutions on iterated products starting from the corresponding
structures on the factors, and give some examples of algebras that
can be described within our theory. We prove a certain result
(called ``invariance under twisting'') for a twisted tensor product
of two algebras, stating that the twisted tensor product does not
change when we apply certain kind of deformation. Under certain
conditions, this invariance can be iterated, containing as
particular cases a number of independent and previously unrelated
results from Hopf algebra theory.
\end{abstract}

\maketitle

\section*{Introduction}
\setcounter{equation}{0}

The difficulty of constructing concrete, nontrivial examples of
noncommutative spaces starting from simpler ones is a common problem
in all different descriptions of noncommutative geometry. If we
think of the commutative situation, we have an easy procedure, the
cartesian product, which allows us to generate spaces of dimension
as big as we want from lower dimensional spaces. Thinking in terms
of the existing dualities between the categories of spaces and the
categories of (commutative) algebras, the natural replacement for
the cartesian product of commutative spaces turns out to be the
tensor product of commutative algebras. The tensor product has often
been considered a replacement for the product of spaces represented
by noncommutative algebras. As it was pointed out in \cite{Cap95a},
this is a very restricted approach. If the ``axiom'' of
noncommutative geometry consists in considering noncommutative
algebras as the representatives for the algebras of functions over
certain ``quantum'' spaces, hence assuming that two different
measurements (or functions) on this kind of spaces do not commute to
each other, then why should we assume that the measurements on the
product commute to each other? There is no reason for imposing this
artificial commutation, hence what we need is a ``noncommutative''
replacement of the tensor product of two algebras, which is supposed
to fit better as an analogue of the product of two noncommutative
spaces and in the same time to be a useful tool for overcoming the
lack of examples formerly mentioned.

When we impose the natural restrictions a product should have,
namely that it contains the factors in a natural way and having
linear size equal to the product of the linear sizes of the factors,
we arrive precisely at the definition of a twisted tensor product
formerly studied by many people, either for the particular case of
algebras (cf. \cite{Tambara90a}, \cite{Cap95a}, \cite{VanDaele94a})
or aiming to define similar structures for discrete groups, Lie
groups, Lie algebras and Hopf algebras (as in \cite{Takeuchi81a},
\cite{Majid90a} and \cite{Michor90a}). Often, this structure appears
in the so-called \emph{factorization problem} of studying under what
conditions we may write an object as a product of two subobjects
having minimal intersection (see for instance the early paper
\cite{Majid90b}) . From a purely algebraic point of view, twisted
tensor products arise as a tool for building algebras starting with
simpler ones, and also, as shown in \cite{VanDaele94a}, in close
relation with certain nonlinear equations. Historically, the
starting point of this theory was the ``braided geometry'' developed
by Majid in the early 1990's, including the ``braided tensor
product'' of algebras in a braided monoidal category, of which the
twisted tensor product of algebras is a sort of ``local'' version.

Whatever the chosen approach to twisted tensor products is, a number
of examples of both classical and recently defined objects fits into
this construction. Ordinary and graded tensor products, crossed
products, Ore extensions and skew group algebras are just some
examples of well-known constructions in classical ring theory that
can be described as twisted tensor products. In the Hopf algebras
and quantum groups area we find smash products, Drinfeld and
Heisenberg doubles, and diagonal crossed products. With a more
geometrical flavour, quantum planes and tori may be realized as
noncommutative products of commutative spaces. And last, but not
least, we may also find some physical models for which this
structure is particularly well suited, such as the Fock space
representations of a particle system with generalized statistics,
which is studied in \cite{Borowiec00a} using techniques which arise
directly from the realisation of certain crossed enveloping algebras
as twisted tensor products.

In the present work, our aim is to look at the twisted tensor
product structure from a more geometrical point of view, regarding
it as the natural representative for the cartesian product of
noncommutative spaces. When we think of this construction
geometrically, it becomes unnatural to restrict ourselves to take
the product of only two spaces, so it appears the problem of finding
suitable conditions that allow us to iterate the construction, and,
whenever this is possible, to check that the obtained iterated
product is ``associative'' in the same sense in which the usual
tensor product is. Also, we will be interested in analyzing whether
we may lift geometrical invariants that we are able to calculate on
the single factors to the iterated twisted product and how to do
this, if possible.

Being such an ubiquitous construction, there are several equivalent
definitions of the twisted tensor product appearing in the
literature, often using different names and notation. In the
Preliminaries we recall some of the results we will use later on,
fixing a unified notation. Concretely, we introduce the definition
of a twisted tensor product $A\otimes_R B$ of two algebras $A$ and
$B$ by means of a twisting map $R:B\otimes A\to A\otimes B$, whose
existence is sufficient for the existence of a deformed product in
the tensor product vector space $A\otimes B$, and is also necessary
when we impose unitality conditions.

In Section \ref{sec:main}, we deal with the problem of iterating the
twisted tensor products, and the lifting of several structures to
the iteration, finding that for three given algebras $A$, $B$ and
$C$, and twisting maps $R_1:B\otimes A\to A\otimes B$, $R_2:C\otimes
B\to B\otimes C$, $R_3:C\otimes A\to A\otimes C$, a sufficient
condition for being able to define twisting maps $T_1:C\otimes
(A\otimes _{R_1}B) \to (A\otimes _{R_1}B)\otimes C$ and
$T_2:(B\otimes _{R_2}C)\otimes A \to A\otimes (B\otimes _{R_2}C)$
associated to $R_1$, $R_2$ and $R_3$ and ensuring that the algebras
$A\otimes_{T_2}(B\otimes_{R_2}C)$ and
$(A\otimes_{R_1}B)\otimes_{T_1}C$ are equal, can be given in terms
of the twisting maps $R_1$, $R_2$ and $R_3$ only. Namely, they have
to satisfy the compatibility condition
\[ (A\otimes R_2)\circ (R_3\otimes B)\circ (C\otimes R_1)=
    (R_1\otimes C)\circ (B\otimes R_3)\circ (R_2\otimes A).
\]
This relation may be regarded as a ``local'' version of the
hexagonal relation satisfied by the braiding of a (strict) braided
monoidal category. We also prove that whenever the algebras and the
twisting maps are unital, the compatibility condition is also
necessary. As it happens for the classical tensor product, and for
the twisted tensor product, the iterated twisted tensor product also
satisfies a Universal Property, which we will state formally in
Theorem \ref{iterunivprop3}. Once the conditions needed to iterate
the construction of the twisted tensor product are fulfilled, we
will prove the Coherence Theorem, stating that whenever one can
build the iterated twisted product of any three factors, it is
possible to construct the iterated twisted product of any number of
factors, and that all the ways one might do this are essentially the
same. This result will allow us to lift to any iterated product
every property that can be lifted to three-factors iterated
products. As applications of the former results  we will
characterize the modules over an iterated twisted tensor product,
also giving a method to build some of them from modules given over
each factor. As a first step towards our aim of building geometrical
invariants over these structures, we will show how to build the
algebras of differential forms and how to lift the involutions of
$\ast$--algebras to the iterated twisted tensor products.

In Section \ref{sec:examples}, we illustrate our theory by
presenting some examples of different structures that arose in
different areas of mathematics and can be constructed using our
method. Two of them (the generalized smash products and diagonal
crossed products) come from Hopf algebra theory, while the other two
(the noncommutative $2n$--planes defined by Connes and
Dubois--Violette, and the observable algebra $\mathcal{A}$ of
Nill--Szlachányi) appear in a more geometrical or physical context.
In particular, we show that the algebras defined by Connes and
Dubois--Violette can be seen as (iterated) noncommutative products
of commutative algebras (as it happens for the quantum planes and
tori), and give a new proof of the fact that the algebra
$\mathcal{A}$ is an AF--algebra, proof which does not imply
calculating any representation of it.  We would like to point out
that the earliest nontrivial example of an iterated twisted tensor
product of algebras was given by Majid in \cite{Majid90c}, in the
form of an iterated sequence of double cross products of certain
bialgebras.

Section \ref{sec:invar} (together with several results from Section
\ref{sec:main}), illustrates the fact that Hopf algebra theory
represents not only a rich source of examples for the theory of
twisted tensor products of algebras, but also a valuable source of
inspiration for it. In this section we prove a result, called
``invariance under twisting'', for a twisted tensor product of two
algebras, which arose as a generalization of the invariance under
twisting for the Hopf smash product (hence the name). It states that
if we start with a twisted tensor product $A\otimes _RB$ together
with a certain kind of datum corresponding to it, we can deform the
multiplication of $A$ to a new algebra structure $A^d$, we can
deform $R$ to a new twisting map $R^d:B\otimes A^d\rightarrow
A^d\otimes B$, so that the twisted tensor products $A^d\otimes
_{R^d}B$ and $A\otimes _RB$ are isomorphic. It turns out that our
result is general enough to include as particular cases some more
independent results from Hopf algebra theory: the well-known theorem
of Majid stating that the Drinfeld double of a quasitriangular Hopf
algebra is isomorphic to an ordinary smash product, a
recent result of Fiore--Steinacker--Wess from \cite{Fiore03a}
concerning a situation where a braided tensor product can be
``unbraided'', and also a recent result of Fiore from \cite{fiore} concerning a situation where a smash product can be ``decoupled'' (this result in turn contains as a particular case the well--known fact that a smash product corresponding to a strongly inner action is isomorphic to the ordinary tensor product).
 We also prove that, under certain circumstances, our
theorem can be iterated, containing thus, as a particular case, the
invariance under twisting of the two-sided smash product from
\cite{BulacuUNa}.

Though we are mainly interested in results of geometrical nature,
and hence most algebras we would like to work with are defined over
the field $\c$ of complex numbers, most of the results can be stated
with no change for algebras over a field or commutative ring $k$,
that we assume fixed throughout all the paper. All algebras will be
supposed to be associative, and usually unital, $k$--algebras. The
term \emph{linear} will always mean \emph{$k$--linear}, and the
unadorned tensor product $\otimes$ will stand for the usual tensor
product over $k$. We will also identify every object with the
identity map defined on it, so that $A\otimes f$ will mean
$\Id_A\otimes f$. For an algebra $A$ we will write $\mu_A$ to denote
the product in $A$ and $u_A:k\rightarrow A$ its unit, and for an
$A$--module $M$ we will use $\l_M$ to denote the action of $A$ on
$M$. For bialgebras and Hopf algebras we use the Sweedler-type
notation $\Delta (h)=h_1\otimes h_2$.

It is worth noting that the proofs of most of our main results are
still valid if instead of considering algebras over $k$ we take
algebras in an arbitrary monoidally closed category.

\section{Preliminaries}\label{sec:prelim}
\setcounter{equation}{0}

\subsection{Twisted tensor products of algebras}

%
%

The notion of twisted tensor product of algebras has been
independently discovered a number of times, and can be found in the
literature under different names and notation. In this section we
collect some results that will be used later, fixing a unified
notation. Main references for definitions and proofs are
\cite{Cap95a} and \cite{VanDaele94a}.

When dealing with spaces that involve a number of tensor products,
notation often becomes obscure and complex. In order to overcome
this difficulty, especially when dealing with iterated products, we
will use a graphical braiding notation in which tangle diagrams
represent morphisms in monoidal categories. For this braiding
notation we refer to \cite{Reshetikhin90a}, \cite{Majid94a} and
\cite{Kassel95a}.

In this notation, a linear map $f:A\to B$ is simply represented by
$\xy 0;/r.125pc/: (0,9)*{\scriptscriptstyle A}; (0,0)*{\DownL{f}};
(0,-8)*{\scriptscriptstyle B};
\endxy$. The composition of morphisms can be written simply by placing
the boxes corresponding to each morphism along the same string,
being the topmost box the corresponding to the map that is applied
in the first place. Several strings placed aside will represent a
tensor product of vector spaces (usually algebras), and a tensor
product of two linear maps, $f\otimes g:A\otimes B\to C\otimes D$
will be written as $\xy 0;/r.125pc/: (0,9)*{\scriptscriptstyle A};
(0,0)*{\DownL{f}}; (0,-8)*{\scriptscriptstyle C};
(8,9)*{\scriptscriptstyle B}; (8,0)*{\DownL{g}};
(8,-8)*{\scriptscriptstyle D};
\endxy$.

With this notation, some well-known properties of morphisms on
tensor products become very intuitive. For instance, the identity
$f\otimes g=(f\otimes D)\circ (A\otimes g)=(C\otimes g)\circ
(f\otimes B)$ is written in braiding notation as
\[ \includegraphics{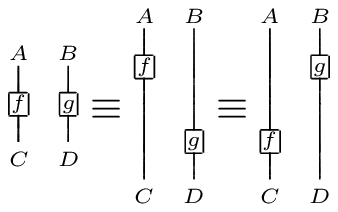}
\]

There are several special classes of morphisms that will receive a
particular treatment. Namely, the identity will be simply written as
a straight line (without any box on it), the algebra product will be
denoted by $\xy 0;/r.125pc/: (0,10)*{\scriptscriptstyle{A}};
(8,10)*{\scriptscriptstyle{A}}; (4,0)*{\Mult};
(4,-8)*{\scriptscriptstyle{A}};
\endxy$. With this notation, the associativity of the algebra
product can be written as:
\[ \includegraphics{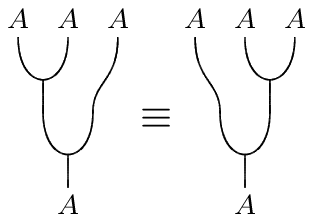}
\]
and the fact that $f:A\to B$ is an algebra morphism may be drawn as
\[ \includegraphics{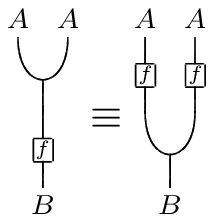}
\]
We will also adopt the convention of not writing the base field (or
ring) whenever it appears as a factor (representing the fact that
scalars can be pushed in or out every factor). According to this
convention, the unit map of an algebra $A$ is represented by $\xy
0;/r.125pc/: (0,-2)*{\Unit{A}}; (0,-8)*{\scriptscriptstyle{A}};
\endxy$, and the compatibility of the unit with the product and with
algebra morphisms are respectively written as
\[ \includegraphics{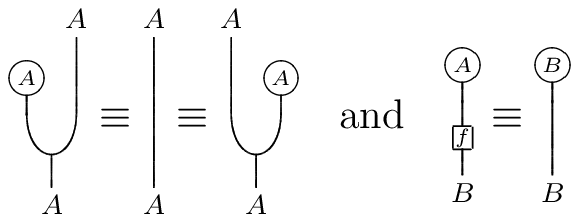}
\]
This conventions may also be applied to module morphisms. If $M$ is
a left $A$--module, we will denote by $\xy 0;/r.125pc/:
(0,10)*{\scriptscriptstyle{A}}; (8,10)*{\scriptscriptstyle{M}};
(4,0)*{\Mult}; (4,-8)*{\scriptscriptstyle{M}};
\endxy$ the module action.
Note that, in spite of the fact that the drawing is the same, there
is no risk of confusing the module action with the algebra product,
since the strings are labeled. Note that, for a morphism $f:M\to N$
of left $A$--modules, the module morphism property is not written
the same way as the algebra morphism property, but as
\[ \includegraphics{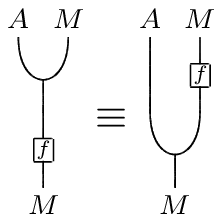}
\]


Recall that given two algebras $A$, $B$ over $k$ and $R:B\otimes
A\to A\otimes B$ a $k$--linear map such that
\begin{eqnarray}
        \label{pentagonequation1} \label{tw2}
        R\circ(B\otimes \mu_A) & = & (\mu_A \otimes B)\circ (A\otimes R)
\circ (R\otimes
        A),\\\label{pentagonequation2} \label{tw3}
        R\circ(\mu_B\otimes A) & = & (A\otimes \mu_B)\circ (R\otimes B)
\circ (B\otimes
        R),
\end{eqnarray}
then the application $\mu_R:=(\mu_A\otimes \mu_B)\circ (A\otimes
R\otimes B)$ is an associative product on $A\otimes B$. In this
case, the map $R$ is said to be a \dtext{twisting map}, and we will
denote by $A\otimes_R B$ the algebra $(A\otimes B,\mu_R)$ that has
$A\otimes B$ as underlying vector space, endowed with the product
$\mu_R$. If, using a Sweedler-type notation, we denote by $R(b\ot
a)=a_R\ot b_R=a_r\otimes b_r$, for $a\in A$, $b\in B$, then
(\ref{tw2}) and (\ref{tw3}) may be rewritten as:
\begin{eqnarray}
&&(aa')_R\ot b_R=a_Ra'_r\ot (b_R)_r, \label{tw4} \\
&&a_R\ot (bb')_R=(a_R)_r\ot b_rb'_R. \label{tw5}
\end{eqnarray}

In braiding notation, we will represent a twisting map $R:B\otimes
A\to A\otimes B$ by a crossing $\xy 0;/r.125pc/:
(0,10)*{\scriptstyle B}; (8,10)*{\scriptstyle A};
(4,0)*{\CruceT{\scriptscriptstyle R}}; (0,-8)*{\scriptstyle A};
(8,-8)*{\scriptstyle B};
\endxy$, where we will omit the label $R$ when there is no risk of confusion, and
equations \eqref{pentagonequation1} and \eqref{pentagonequation2}
are represented respectively by
\[ \includegraphics{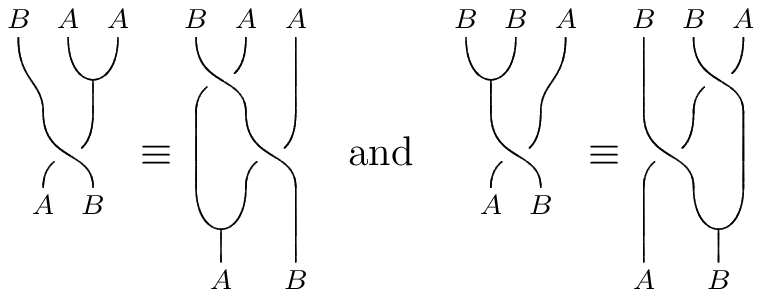}
\]

For further use, we record the following consequence of (\ref{tw4})
and (\ref{tw5}):
\begin{eqnarray}
(aa')_R\ot (bb')_R=(a_R)_{\mathcal{R}}(a'_r)_{\overline{r}}\ot
(b_{\mathcal{R}})_{\overline{r}}(b'_R)_r, \label{patru}
\end{eqnarray}
for all $a, a'\in A$ and $b, b'\in B$, where $\mathcal{R}$ and
$\overline{r}$ are two more copies of $R$; in braiding notation this
last identity is written as:
\[ \includegraphics{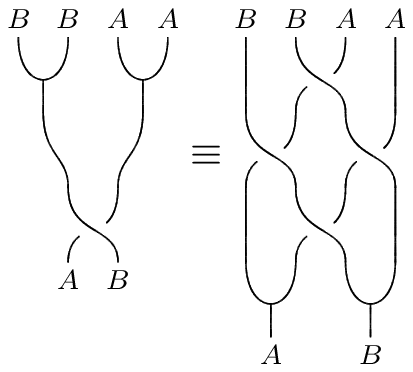}
\]

Whenever $A$ and $B$ are unital, if $R$ is a twisting map that
satisfies the extra conditions
\begin{equation}\label{unitconditions}
\left.
    \begin{array}{rcl}
        R(1\otimes a)&=&a\otimes 1\\
        R(b\otimes 1)&=&1\otimes b
\end{array}\right\}
\end{equation}
then the canonical maps $i_A:A\to A\otimes_R B$ and $i_B:B\to
A\otimes_R B$ defined by $i_A(a):=a\otimes 1$, $i_B(b):=1\otimes b$,
are algebra morphisms, and $A\otimes_R B$ is a unital algebra, with
unit $1\otimes 1$. In this case, we say that $R$ is a \dtext{unital
twisting map}. Most of the twisting maps we will study are unital;
however, it is worth noting that associativity constraints do not
depend on the unitality of the twisting map. In braiding notation,
the unitality conditions read
\[ \includegraphics{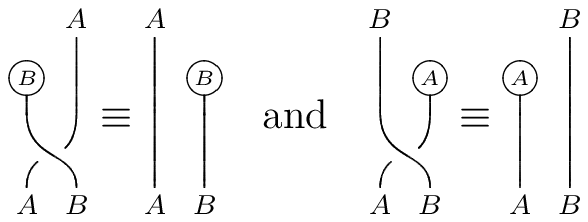}
\]

A special family of examples of twisting maps involves bijective
maps. Concerning this situation, we can state the following result from \cite{cmz},
which will be used later:

\begin{proposition} \label{bij}
Let $A\ot _RB$ be a twisted tensor product of algebras such that the
map $R$ is bijective, and denote by $V:A\ot B\rightarrow B\ot A$ its
inverse. Then $V$ is also a twisting map and $R$ is an algebra
isomorphism between $B\ot_V A$ and $A\ot_R B$.
\end{proposition}


In classical homological algebra, the usual tensor product is
commonly introduced by means of its universal property, where the
commutation between elements belonging to the first factor and
elements belonging to the second one is implicitly required. In this
property, we have to consider the canonical algebra monomorphisms
$i_A:A\hookrightarrow A\otimes B$ and $i_B:B\hookrightarrow A\otimes
B$ given by $i_A(a):=a\otimes 1$ and $i_B(b):=1\otimes b$
respectively. Because of the twisting map conditions, these maps are
still algebra morphisms when we consider a twisted tensor product
$A\otimes_R B$ instead of $A\otimes B$; moreover, twisted tensor
products may be characterized as algebra structures defined on
$A\otimes B$ such that the above maps are algebra inclusions and
satisfying $a\otimes b= i_A(a)i_B(b)$ for all $a\in A$, $b\in B$. As
a consequence, with a slight modification, that essentially involves
replacing the usual flip by the twisting map, one may also state a
universal property for twisted tensor products, as shown in
\cite{Caenepeel00a}:

\begin{theorem}\label{twistuniversalproperty}
Let $A$, $B$ be two $k$--algebras, and let $R:B\otimes A\to A\otimes
B$ be a unital twisting map. Given a $k$--algebra $X$, and algebra
morphisms $u:A\to X$, $v:B\to X$ such that
    \begin{equation}\label{univpropertycondition}
        \mu_X\circ (v\otimes u) = \mu_X\circ (u\otimes v)\circ R,
    \end{equation}
then we can find a unique algebra map $\vphi:A\otimes_R B\to X$ such
that
    \begin{align}
        \vphi\circ i_A & =  u, \\
        \vphi \circ i_B &= v.
    \end{align}
\end{theorem}

If $A$ and $B$ are $\ast$--algebras with involutions $j_A$ and
$j_B$, and $R:B\otimes A\to A\otimes B$ is a twisting map such that
\begin{equation}\label{starcondition}
    (R\circ(j_B\otimes j_A)\circ\tau) \circ (R\circ(j_B\otimes
    j_A)\circ\tau)=A\otimes B,
\end{equation}
then $A\otimes_R B$ is a $\ast$--algebra with involution
$R\circ(j_B\otimes j_A)\circ\tau$, where $\tau:A\otimes B\to
B\otimes A$ denotes the usual flip. Moreover, if $R$ is unital, then
$i_A$ and $i_B$ become $\ast$--morphisms. This involutive condition
is written down in braiding notation in the following way:
\[ \includegraphics{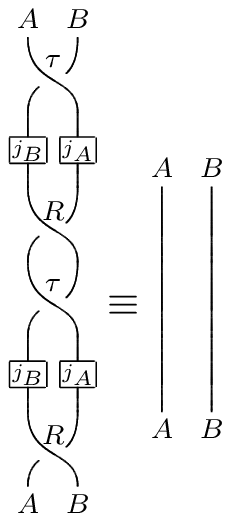}
\]

When we have a left $A$--module $M$, a left $B$--module $N$, a
twisting map $R:B\otimes A\to A\otimes B$ and a linear map
$\tau_{M,B}:B\otimes M \to M\otimes B$ such that
\begin{eqnarray}\label{moduletwistcondition1}
        \tau_{M,B}\circ(\mu_B\otimes M) & = & (M\otimes \mu_B)\circ (\tau_{M,B}\otimes B)\circ (B\otimes
        \tau_{M,B}), \\ \label{moduletwistcondition2}
        \tau_{M,B}\circ(B\otimes \l_M) & = & (\l_M \otimes B)\circ (A\otimes \tau_{M,B})\circ (R\otimes
        M),
\end{eqnarray}
then the map $\l_{\tau_{M,B}}:(A\otimes_R B)\otimes (M\otimes N)\to
M\otimes N$ defined by $\l_{\tau_{M,B}}:=(\l_M\otimes \l_N)\circ
(A\otimes \tau_{M,B}\otimes N)$ yields a left $(A\otimes_R
B)$--module structure on $M\otimes N$, which furthermore is
compatible with the inclusion of $A$. In this case, we say that
$\tau_{M,B}$ is a \dtext{(left) module twisting map}. If we denote
by $\xy 0;/r.125pc/: (0,10)*{\scriptstyle B}; (8,10)*{\scriptstyle
M}; (4,0)*{\Cruce}; (0,-8)*{\scriptstyle M}; (8,-8)*{\scriptstyle
B};
\endxy$ the module twisting map, the module twisting conditions
look the same as the twisting conditions for algebra twisting maps
(replacing $A$ by $M$). Unlike what happens for algebra twisting
maps, usually is not enough to have a left $(A\otimes_R B)$--module
structure on $M\otimes N$ in order to recover a module twisting map.
Some sufficient conditions for this to happen can be found in
\cite{Cap95a}

Besides module lifting conditions, in \cite{Cap95a} is shown how to
lift twisting maps to algebras of differential forms on them. More
precisely:
\begin{theorem}\label{twistdifferentialforms}
Let $A$, $B$ be two algebras. Then any twisting map $R:B\otimes A\to
A\otimes B$ extends to a unique twisting map $\tilde{R}:\O B\otimes
\O A \to \O A\otimes \O B$ which satisfies the conditions
\begin{eqnarray}\label{twisteddiff1}
\tilde{R}\circ(d_B\otimes \O A)&=&(\eps_A \otimes d_B)\circ
\tilde{R}, \\
\label{twisteddiff2} \tilde{R}\circ (\O B\otimes d_A)&=&(d_A \otimes
\eps_B)\circ \tilde{R},
\end{eqnarray}
where $d_A$ and $d_B$ denote the differentials on $\O A$ and $\O B$,
and $\eps_A$, $\eps_B$ stand for the gradings on $\O A$ and $\O B$,
respectively. Moreover, $\O A\otimes_{\tilde{R}}\O B$ is a graded
differential algebra with differential $d(\vphi\otimes
\omega):=d_A\vphi\otimes \omega + (-1)^{\abs{\vphi}}\vphi\otimes
d_B\omega$.
\end{theorem}

Conditions \eqref{twisteddiff1} and \eqref{twisteddiff2} can be
translated, in braiding notation, to the equalities
\[ \includegraphics{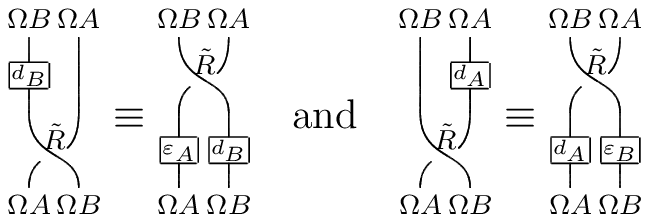}
\]
respectively.

\subsection{The noncommutative planes of Connes and Dubois--Violette}

The original definition of noncommutative $4$--planes (and
$3$--spheres) arises from some $K$--theoretic equations, inspired by
the properties of the Bott projector on the cohomology of classical
spheres. We do not need this interpretation here, so we adopt
directly the equivalent definition given by means of generators and
relations. Any reader interested in full details on the construction
and properties of noncommutative planes and spheres should look at
\cite{Connes02a}. Our study will be centered on the noncommutative
planes associated to critical points of the scaling foliation, as
the definition of the noncommutative plane in these points is easily
generalized to higher dimensional frameworks.

Let us then consider $\theta\in \M_n(\R)$ an antisymmetric matrix,
$\theta=(\theta_{\mu\nu})$, $\theta_{\nu\mu}=-\theta_{\mu\nu}$, and
let $C_{alg}(\R^{2n}_\theta)$ be the associative algebra generated
by $2n$ elements $\{z^\mu,\bar{z}^\mu\}_{\mu=1,\dotsc,n}$ with
relations
\begin{equation} \label{ncplanerelations}
\left.\begin{array}{r}
z^\mu z^\nu=\l^{\mu\nu}z^\nu z^\mu \\
\bar{z}^\mu \bar{z}^\nu=\l^{\mu\nu}\bar{z}^\nu \bar{z}^\mu \\
\bar{z}^\mu z^\nu = \l^{\nu \mu} z^\nu \bar{z}^\mu
\end{array}\right\}\forall\,\mu,\nu=1,\dotsc,n,\ \text{being
$\l^{\mu\nu}:=e^{i\theta_{\mu\nu}}$.}
\end{equation}
Note that $\l^{\nu\mu}=(\l^{\mu\nu})^{-1}=\overline{\l^{\mu\nu}}$
for $\mu\neq \nu$, and $\l^{\mu\mu}=1$ by antisymmetry.

We can now endow the algebra $C_{alg}(\R^{2n}_\theta)$ with the
unique involution of $\c$--algebras $x\mapsto x^\ast$ such that
$(z^\mu)^\ast=\bar{z}^\mu$. This involution gives a structure of
$\ast$--algebra on $C_{alg}(\R^{2n}_\theta)$. As a $\ast$--algebra,
$C_{alg}(\R^{2n}_\theta)$ is a deformation of the commutative
algebra $C_{alg}(\R^{2n})$ of complex polynomial functions on
$\R^{2n}$, and it reduces to it when we take $\theta=0$. The algebra
$C_{alg}(\R^{2n}_\theta)$ will be then referred to as the
\dtext{(algebra of complex polynomial functions on the)
noncommutative $2n$--plane $\R^{2n}_\theta$}. In fact, former
relations define a deformation $\c^n_\theta$ of $\c^n$, so we can
identify the noncommutative complex $n$--plane $\c^n_\theta$ with
$\R^{2n}_\theta$ by writing
$C_{alg}(\c^n_\theta):=C_{alg}(\R^{2n}_\theta)$.

We define $\O_{alg}(\R^{2n}_\theta)$, \dtext{the algebra of
algebraic differential forms on the noncommutative plane
$\R^{2n}_\theta$}\index{differential forms}, to be the complex
unital associative graded algebra
\[ \O_{alg}(\R^{2n}_\theta):=\bigoplus_{p\in\n}\O_{alg}^p(\R^{2n}_\theta)
\]
generated by $2n$ elements $z^\mu$, $\bar{z}^\mu$ of degree 0, with
relations:
\[
\left.\begin{array}{r}
z^\mu z^\nu=\l^{\mu\nu}z^\nu z^\mu \\
\bar{z}^\mu \bar{z}^\nu=\l^{\mu\nu}\bar{z}^\nu \bar{z}^\mu \\
\bar{z}^\mu z^\nu = \l^{\nu \mu} z^\nu \bar{z}^\mu
\end{array}\right\}\forall\,\mu,\nu=1,\dotsc,n,\ \text{being
$\l^{\mu\nu}:=e^{i\theta_{\mu\nu}}$,}
\]
and by $2n$ elements $dz^\mu$, $d\bar{z}^\mu$ of degree 1, with
relations:
\begin{equation}
\label{relations6.1-2} \left.
\begin{array}{rcl}
dz^\mu dz^\nu+\l^{\mu\nu}dz^\nu dz^\mu  &= & 0, \\
d\bar{z}^\mu d\bar{z}^\nu+\l^{\mu\nu}d\bar{z}^\nu d\bar{z}^\mu  &= & 0,\\
d\bar{z}^\mu dz^\nu+\l^{\nu\mu}dz^\nu d\bar{z}^\mu  & = & 0,\\
\end{array}\
\begin{array}{rcl}
z^\mu dz^\nu & = & \l^{\mu\nu}dz^\nu z^\mu, \\
\bar{z}^\mu d\bar{z}^\nu & = & \l^{\mu\nu}d\bar{z}^\nu \bar{z}^\mu,\\
\bar{z}^\mu dz^\nu & = & \l^{\nu\mu}dz^\nu \bar{z}^\mu, \\
z^\mu d\bar{z}^\nu & = & \l^{\nu\mu}d\bar{z}^\nu z^\mu,
\end{array} \right\}\
\text{$\forall\;\mu,\nu = 1,\dotsc, n$}.
\end{equation}

In this setting, there exists a unique differential $d$ of
$\O_{alg}(\R^{2n}_\theta)$ (that is, an antiderivation of degree 1
such that $d^2=0$) which extends the mapping $z^\mu\mapsto dz^\mu$,
$\bar{z}^\mu\mapsto d\bar{z}^\mu$. Indeed, such a differential is
obtained by extending the definition on the generators according to
the Leibniz rule. With this differential, $\O_{alg}(\R^{2n}_\theta)$
becomes a graded differential algebra. It is also possible to extend
the mapping $z^\mu\mapsto \bar{z}^\mu$, $dz^\mu\mapsto
d\bar{z}^\mu=:\overline{(dz^\mu)}$ to the whole algebra
$\O_{alg}(\R^{2n}_\theta)$ as an antilinear involution
$\omega\mapsto \overline{\omega}$ such that
$\overline{\omega\omega'}=(-1)^{pq}\overline{\omega'}\overline{\omega}$
for any $\omega\in \O_{alg}^p(\R^{2n}_\theta)$, $\omega'\in
\O_{alg}^q(\R^{2n}_\theta)$. For this extension we have that
$d\overline{\omega}=\overline{d\omega}$.

Our interest in these algebras arises from the fact that the
noncommutative 4--plane can easily be realized as a twisted tensor
product of two commutative algebras (namely as a twisted product of
two copies of $\c[x,\bar{x}]$, which is nothing but the algebra of
polynomial functions on the complex plane), hence looking like the
algebra representing a sort of \emph{noncommutative cartesian
product} of two commutative spaces. Our original interest in
iterated twisted tensor products came when we asked ourselves about
the possibility of looking at the $2n$--noncommutative plane as a
certain product of commutative algebras.


\section{Iterated twisted tensor products}\label{sec:main}
\setcounter{equation}{0}

In this section, our aim is to study the construction of iterated
twisted tensor products. If we think of twisted tensor products as
natural noncommutative analogues for the usual cartesian product of
spaces, it is natural to require that the product of three or more
spaces still respects every single factor.

Morally, the construction of a twisting map boils down to giving a
rule for exchanging factors between the algebras involved in the
product. A natural way for doing this would be to perform a series
of two factors twists, that should be related to the already given
notion of twisting map, and afterwards to apply algebra
multiplication in each factor.

Suppose that $A$, $B$ and $C$ are algebras, let
\begin{gather*}
R_1:B\otimes A\lto A\otimes B,\\
R_2:C\otimes B\lto B\otimes C,\\
R_3:C\otimes A\lto A\otimes C
\end{gather*}
(unital) twisting maps, and consider now the application
\[
T_1:C\otimes (A\otimes _{R_1}B) \lto (A\otimes _{R_1}B)\otimes C
\]
given by $T_1:=(A\otimes R_2)\circ (R_3\otimes B)$. We can also
build the map
\[
T_2:(B\otimes _{R_2}C)\otimes A \lto A\otimes (B\otimes _{R_2}C)
\]
given by $T_2=(R_1\otimes C)\circ (B\otimes R_3)$. It is a natural
question to ask if these maps are twisting maps. In general, this is
not the case, as we will show in (Counter)example
\ref{florincounterexample1}. In the following Theorem, we state
necessary and sufficient conditions for this to happen.

\begin{theorem}\label{itertwisting}
With the above notation, the following conditions are equivalent:
\begin{enumerate}
    \item $T_1$ is a twisting map.
    \item $T_2$ is a twisting map.
    \item The maps $R_1$, $R_2$ and $R_3$ satisfy the following compatibility
    condition (called the hexagon equation):
    \begin{equation}\label{hexagonequation} (A\otimes R_2)\circ (R_3\otimes B)\circ (C\otimes R_1)=
    (R_1\otimes C)\circ (B\otimes R_3)\circ (R_2\otimes A),
    \end{equation}
    that is, the following diagram is commutative.
    \[ \xymatrix{
            & C\otimes A\otimes B \ar[r]^{R_3\otimes B}& A\otimes
            C\otimes B \ar[dr]^{A\otimes R_2} \\
            C\otimes B\otimes A \ar[ur]^{C\otimes R_1} \ar[dr]_{R_2\otimes A}
            & & & A\otimes B\otimes C\\
            & B\otimes C\otimes A \ar[r]^{B\otimes R_3}& B\otimes
            A\otimes C \ar[ur]^{R_1\otimes C}
            }
    \]
\end{enumerate}
Moreover, if all the three conditions are satisfied, then the
algebras $A\otimes_{T_2}(B\otimes_{R_2}C)$ and
$(A\otimes_{R_1}B)\otimes_{T_1}C$ are equal. In this case, we will
denote this algebra by $A\otimes_{R_1}B\otimes_{R_2} C$.
\end{theorem}

\begin{proof}
We prove only the equivalence between $(1)$ and $(3)$, being the
equivalence between $(2)$ and $(3)$ completely analogous.

\imp{3}{1} Suppose that the hexagon equation is satisfied. In order
to prove that $T_1$ is a twisting map, we have to check the
conditions \eqref{pentagonequation1} and \eqref{pentagonequation2}
for $T_1$, namely, we have to check the relations
\begin{eqnarray}
T_1\circ(C\otimes \mu_{R_1}) &=& (\mu_{R_1}\otimes C)\circ(A\otimes
B\otimes T_1)\circ(T_1\otimes A\otimes B), \label{twistT1:1}\\
T_1\circ(\mu_C\otimes A\otimes B) &=& (A\otimes B\otimes \mu_C)\circ
(T_1\otimes C)\circ (C\otimes T_1).\label{twistT1:2}
\end{eqnarray}

To prove this we use braiding notation. Taking into account that the
hexagon equation is written as:
\[ \includegraphics{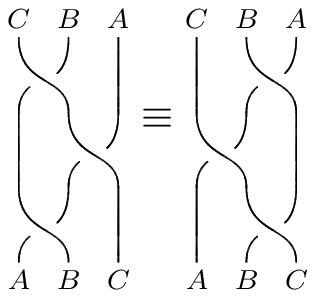}
\]
the proof of condition \eqref{twistT1:1} is given by:
\[
\includegraphics[scale=1]{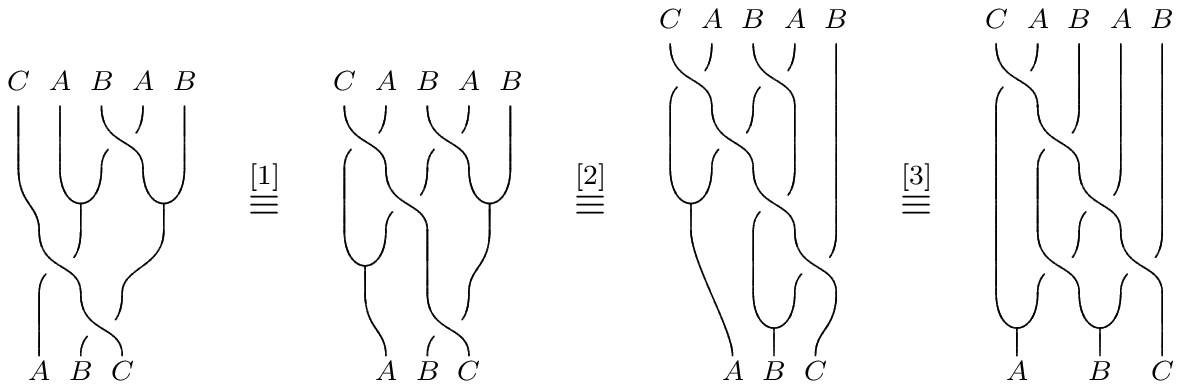}
\]
where in $[1]$ we use the twisting condition for $R_3$, in $[2]$ we
use the twisting condition for $R_2$, and in $[3]$ we use the
hexagon equation. On the other hand, condition \eqref{twistT1:2} is
proven as follows:
\[
\includegraphics[scale=1]{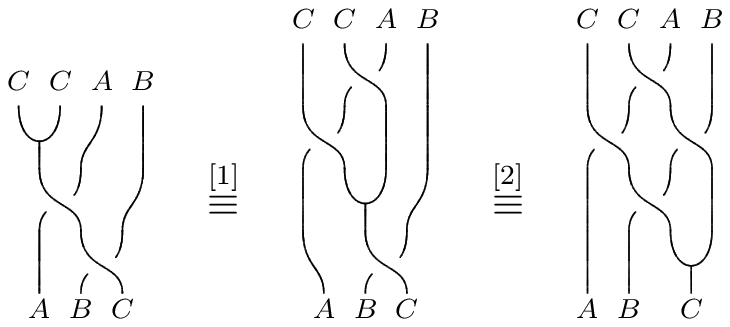}
\]
where now $[1]$ is due to the twisting conditions for $R_3$, and
$[2]$ to twisting conditions for $R_2$. This proves that $T_1$
satisfies the pentagonal equations. Furthermore, if $R_2$ and $R_3$
are unital, then we have that
\begin{gather*}
T_1(c\otimes 1\otimes 1)=(A\otimes R_2)(R_3\otimes B)(c\otimes
1\otimes 1)=(A\otimes R_2)(1\otimes c\otimes 1)= 1\otimes 1\otimes
c,\\
T_1(1\otimes a\otimes b)=(A\otimes R_2)(R_3\otimes B)(1\otimes
a\otimes b)=(A\otimes R_2)(a\otimes 1\otimes b)= a\otimes b\otimes
1,
\end{gather*}
so $T_1$ is also a unital twisting map.

\imp{1}{3} Now we assume \eqref{twistT1:1} and \eqref{twistT1:2}. It
is enough to apply \eqref{twistT1:1} to an element of the form
$c\otimes 1\otimes b\otimes a\otimes 1$ in order to recover the
hexagon equation for a generic element $c\otimes b\otimes a$ of the
tensor product $C\otimes B\otimes A$.

To finish the proof, assume that the three equivalent conditions are
satisfied. To see that the algebras
$A\otimes_{T_2}(B\otimes_{R_2}C)$ and
$(A\otimes_{R_1}B)\otimes_{T_1}C$ are equal, it is enough to expand
the expressions of the products
\begin{eqnarray*}
\mu_{T_2}&=&(\mu_A\otimes \mu_{R_2})\circ(A\otimes T_2\otimes
B\otimes
C), \\
\mu_{T_1}&=&(\mu_{R_1}\otimes \mu_C)\circ (A\otimes B\otimes
T_1\otimes C),
\end{eqnarray*}
and realize that they are exactly the same application, for which we
only have to observe that $ (A\otimes B\otimes R_2)\circ (R_1\otimes
C\otimes B)=R_1\otimes R_2=(R_1\otimes B\otimes C)\circ (B\otimes
A\otimes R_2)$. \qed\end{proof}

When three twisting maps satisfy the hypotheses of Theorem
\ref{itertwisting}, we will say either that they are
\dtext{compatible twisting maps}, or that the twisting maps satisfy
the \dtext{hexagon (or braid) equation}. If the twisting maps $R_i$
are not unital, the hexagon equation is still sufficient for getting
associative products associated to $T_1$ and $T_2$, but in general
we need unitality to recover the compatibility condition from the
associativity of the iterated products.

One could wonder whether the braid relation is automatically
satisfied for any three unital twisting maps. This is not the case,
as shown in the following example:

\begin{example}\label{florincounterexample1}
Take $H$ a \emph{noncocommutative} (finite dimensional) bialgebra,
$A=B=H^\ast$, $C=H$. Consider the left regular action of $H$ on
$H^\ast$ given by $(h\rightharpoonup p)(h'):=p(h'h)$; with this
action, $H^\ast$ becomes a left $H$--module algebra, so we can
define the twisting map induced by the action as:
\begin{eqnarray*}
\sigma:H\otimes H^\ast & \lto & H^\ast\otimes H \\
h\otimes p & \lmto & (h_{1}\rightharpoonup p)\otimes h_{2}.
\end{eqnarray*}
If we consider now the twisting maps $R_1:B\otimes A\lto A\otimes
B$, $R_2:C\otimes B\lto B\otimes C$, $R_3:C\otimes A\lto A\otimes
C$, defined as $R_1:=\tau$, $R_2=R_3:=\sigma$, being $\tau$ the
usual flip, then the braid relation among $R_1$, $R_2$ and $R_3$
boils down to the equality
\[ (h_{1}\rightharpoonup q)\otimes (h_{2}\rightharpoonup p)\otimes h_{3}=
(h_{2}\rightharpoonup q)\otimes (h_{1}\rightharpoonup p)\otimes
h_{3},
\]
for all $h\in H$, $p,q\in H^\ast$, but this relation is false, as we
chose $H$ to be noncocommutative.
\end{example}

\begin{remark} The multiplication in the algebra $A\ot _{R_1}B\ot _{R_2}C$
can be given, using the Sweedler-type notation recalled before, by
the formula:
\begin{eqnarray}
&&(a\ot b \ot c)(a'\ot b'\ot c')=a(a'_{R_3})_{R_1}\ot
b_{R_1}b'_{R_2}\ot (c_{R_3})_{R_2}c'. \label{trei}
\end{eqnarray}
\end{remark}

The next natural question that arises is whether whenever we have a
twisting map $T:C\otimes (A\otimes_R B)\to (A\otimes_R B)\otimes C$,
it splits as a composition of two suitable twisting maps. Once
again, this is not possible in general.

\begin{theorem}[Right splitting]
Let $A$, $B$, $C$ be algebras, $R_1:B\otimes A\to A\otimes B$ and
$T:C\otimes (A\otimes _{R_1}B)\to (A\otimes _{R_1}B) \otimes C$
unital twisting maps. The following are equivalent:
\begin{enumerate}
    \item There exist $R_2:C\otimes B\to B\otimes C$ and
    $R_3:C\otimes A\to A\otimes C$ twisting maps such that
    $T=(A\otimes R_2)\circ (R_3\otimes B)$.
    \item The map $T$ satisfies the \dtext{(right) splitting conditions}:
        \begin{eqnarray}
        T(C\otimes (A\otimes 1)) & \subseteq & (A\otimes 1)\otimes C,\\
        T(C\otimes (1\otimes B)) & \subseteq & (1\otimes B)\otimes C.
        \end{eqnarray}
\end{enumerate}
\end{theorem}

\begin{proof}\hfill\\
\imp{1}{2} It is trivial.

\imp{2}{1} Because of the conditions imposed to $T$, the map
$R_2:C\otimes B\to B\otimes C$ given as the only $k$--linear map
such that $(u_A\otimes R_2)\circ (\tau \otimes B)= T\circ(C\otimes
(u_A\otimes B))$ is well defined. From the fact that $T$ is a
twisting map it is immediately deduced that also $R_2$ is a twisting
map. Analogously, we can define $R_3:C\otimes A\to A\otimes C$ as
the only $k$--linear map such that $(A\otimes \tau)\circ(R_3\otimes
u_B)=T\circ (C\otimes (A\otimes u_B))$, which is also a well defined
twisting map. We only have to check that $T=(A\otimes R_2)\circ
(R_3\otimes B)$. Using braiding notation we have
\[
\includegraphics[scale=1]{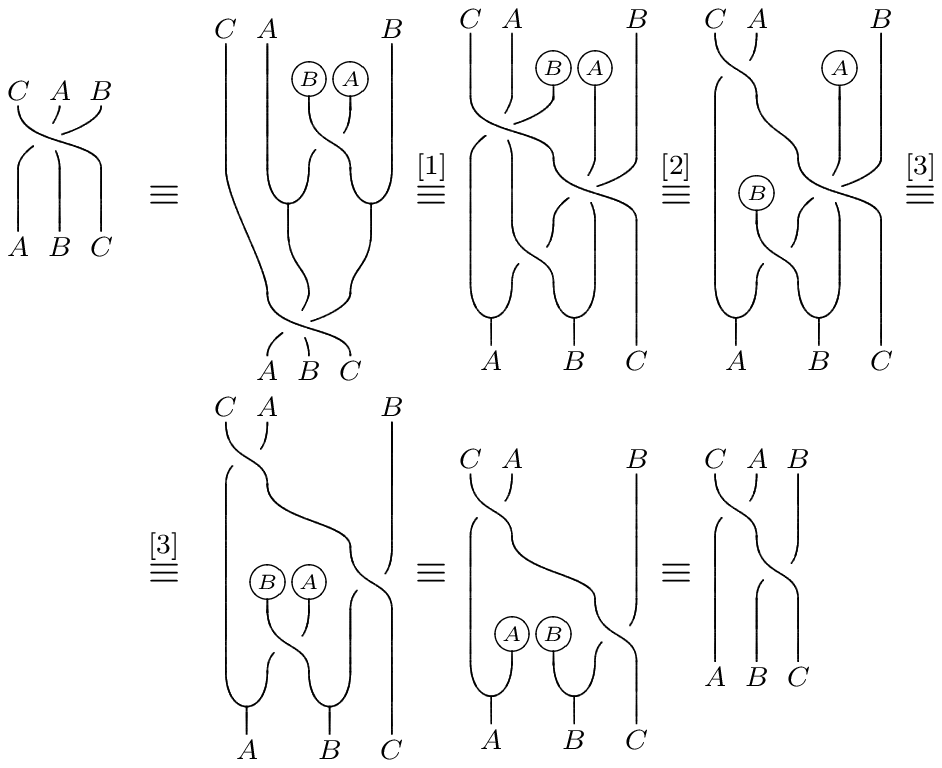}
\]
as we wanted to show, and where in $[1]$ we are using that $T$ is a
twisting map, and in $[2]$ and $[3]$ the definitions of $R_3$ and
$R_2$ respectively. \qed\end{proof}

Again, we can ask ourselves whether the condition we required for
the twisting map $T$ to split might be trivial. The following
example shows a situation in which an iterated twisted tensor
product cannot be split:

\begin{example}
We give an example of twisting maps $R:B\ot A\rightarrow A\ot B$ and
$T:C\ot (A\ot _R B)\rightarrow (A\ot _R B)\ot C$ for which it is
{\it not}
true that $T(c\ot (a\ot 1))\in (A\ot 1)\ot C$ for all $a\in A$, $c\in C$.\\
Let $H$ be a finite dimensional Hopf algebra with antipode $S$.
Recall that the Drinfeld double $D(H)$ is a Hopf algebra having
$H^{* cop}\ot H$ as coalgebra structure and multiplication
\begin{eqnarray*}
&&(p\ot h)(p'\ot h')=p(h_1\rightharpoonup p'\leftharpoonup
S^{-1}(h_3))\ot h_2h',
\end{eqnarray*}
for all $p, p'\in H^*$ and $h, h'\in H$, where $\rightharpoonup $
and $\leftharpoonup $ are the left and right regular actions of $H$
on $H^*$ given by $(h\rightharpoonup p)(h')=p(h'h)$ and
$(p\leftharpoonup h)(h')= p(hh')$. The Heisenberg double
${\mathcal{H}}(H)$ is the smash product $H\# H^*$, where $H^*$ acts
on $H$ via the left regular action $p\rightharpoonup h= p(h_2)h_1$.
Recall from \cite{Lu94a} that ${\mathcal{H}}(H)$ becomes a left
$D(H)$-module algebra, with action
\begin{eqnarray*}
&&(p\ot h)\rightharpoonup (h'\ot q)=p_2(h'_2)q_2(h)(h'_1\ot p_3q_1
S^{* -1}(p_1)),
\end{eqnarray*}
for all $p, q\in H^*$ and $h, h'\in H$, which is just the left
regular action of $D(H)$ on ${\mathcal{H}}(H)$ identified as vector
space with
$D(H)^*$. \\
Now, we take $A=H$, $B=H^*$, $C=D(H)$, $R:H^*\ot H\rightarrow H\ot
H^*$, $R(p\ot h)=p_1\rightharpoonup h\ot p_2$ (hence $H\ot _R
H^*=H\# H^*= {\mathcal{H}}(H)$) and
\begin{eqnarray*}
&&T:D(H)\ot {\mathcal{H}}(H)\rightarrow {\mathcal{H}}(H)\ot D(H), \\
&&T((p\ot h)\ot (h'\ot q))=(p\ot h)_1\rightharpoonup (h'\ot q) \ot
(p\ot h)_2
\end{eqnarray*}
(hence ${\mathcal{H}}(H)\ot _TD(H)={\mathcal{H}}(H)\# D(H)$, so $T$
is a twisting map). Now we can see that
\begin{eqnarray*}
&&T((p\ot h)\ot (h'\ot 1))=p_3(h'_2)(h'_1\ot p_4S^{* -1}(p_2))\ot
(p_1\ot h),
\end{eqnarray*}
which in general does {\it not} belong to $(H\ot 1)\ot D(H)$.
\end{example}

Of course, there exists an analogous left splitting theorem, that we
state for completeness, and whose proof is analogous to the former
one.

\begin{theorem}[Left splitting]
Let $A$, $B$, $C$ be algebras, $R_2:C\otimes B\to B\otimes C$ and
$T:(B\otimes _{R_2}C)\otimes A \to A\otimes (B \otimes _{R_2}C)$
twisting maps. The following are equivalent:
\begin{enumerate}
    \item There exist $R_1:B\otimes A\to A\otimes B$ and
    $R_3:C\otimes A\to A\otimes C$ twisting maps such that
    $T=(R_1\otimes C)\circ (B\otimes R_3)$.
    \item The map $T$ satisfies the \dtext{(left) splitting conditions}:
        \begin{eqnarray}
            T((1\otimes C)\otimes A) & \subseteq & (A\otimes 1)\otimes C,\\
            T((B\otimes 1)\otimes A) & \subseteq & A\otimes (B\otimes 1).
        \end{eqnarray}
\end{enumerate}
\end{theorem}

The universal property (Theorem \ref{twistuniversalproperty})
formerly stated can be easily extended to the iterated setting, as
we show in the following result:

%

\begin{theorem} \label{iterunivprop3} Let $(A,B,C,R_1,R_2,R_3)$ be as in Theorem
\ref{itertwisting}. Assume that we have a $k$--algebra $X$ and
algebra morphisms $u:A \to X$, $v:B \to X$, $w:C \to X$, such that
\begin{equation}\label{iterunivpropcondition}
    \mu_X\circ(w\otimes v\otimes u) = \mu_X\circ (u\otimes v\otimes
    w)\circ (A\otimes R_2)\circ (R_3\otimes B)\circ (C\otimes R_1).
\end{equation}
Then there exists a unique algebra map $\vphi:A\otimes_{R_1}
B\otimes_{R_2} C\to X$ such that $\vphi\circ i_A=u$, $\vphi\circ
i_B=v$, $\vphi\circ i_C=w$.
\end{theorem}

\begin{proof}
Assume the we have a map $\vphi$ satisfying the conditions in the
theorem, then we may write
    \begin{eqnarray*}
        \vphi(a\otimes b\otimes c) &=& \vphi((a\otimes 1\otimes 1)(1\otimes b\otimes 1)(1\otimes 1\otimes
        c)) \\
        &=& \vphi(a\otimes 1\otimes 1)\vphi(1\otimes b\otimes 1)\vphi(1\otimes 1\otimes
        c)) \\
        &=& \vphi (i_A(a))\vphi(i_B(b))\vphi(i_C(c))
        \\
        &=& u(a)v(b)w(c),
    \end{eqnarray*}
and so $\vphi$ is uniquely defined.

For the existence part, define $\vphi(a\otimes b\otimes
c):=u(a)v(b)w(c)$, and let us check that this map is indeed an
algebra morphism. Using formula \eqref{trei}, we have
    \begin{eqnarray*}
        \vphi((a\otimes b\otimes c)(a'\otimes b'\otimes c')) & =
        & \vphi (a(a'_{R_3})_{R_1}\ot b_{R_1}b'_{R_2}\ot (c_{R_3})_{R_2}c')\\
        &=& u(a)u((a'_{R_3})_{R_1})v(b_{R_1})v(b'_{R_2})w((c_{R_3})_{R_2})w(c').
    \end{eqnarray*}
On the other hand, we have
    \begin{eqnarray*}
        \vphi(a\ot b\ot c)\vphi(a'\ot b'\ot c') &=&
        u(a)v(b)w(c)u(a')v(b')w(c') \\
        &=& u(a)v(b)u(a'_{R_3})w(c_{R_3})v(b')w(c')\\
        &=&
        u(a)u((a'_{R_3})_{R_1})v(b_{R_1})v(b'_{R_2})w((c_{R_3})_{R_2})w(c'),
    \end{eqnarray*}
and thus we conclude that $\vphi$ is an algebra morphism. The fact
that $\vphi$ satisfies the required relations with $u$, $v$ and $w$
is immediately deduced from its definition. \qed\end{proof}

To reach completely the aim of defining an analogue for the product
of spaces, one should be able to construct a product of any number
of factors. In order to construct the three--factors product, we had
to add one extra condition, namely the hexagon equation, to the
conditions that were imposed for building the two--factors product
(the twisting map conditions). Fortunately, in order to build a
general $n$--factors twisted product of algebras one needs no more
conditions besides the ones we have already met. Morally, this just
means than having pentagonal (twisting) and hexagonal (braiding)
conditions, we can build any product without worrying about where to
put the parentheses. The way to prove this is using induction. As
our induction hypothesis, we assume that whenever we have $n-1$
algebras $B_1,\dotsc, B_{n-1}$, with a twisting map
$S_{ij}:B_j\otimes B_i\to B_i\otimes B_j$ for every $i<j$, and such
that for any $i<j<k$ the maps $S_{ij}$, $S_{jk}$ and $S_{ik}$ are
compatible, then we can build the iterated product
$B_1\otimes_{S_{12}}B_2\otimes_{S_{23}} \dotsb \otimes_{S_{n-1\, n}}
B_n$ without worrying about parentheses. Let then $A_1,\dotsc,A_n$
be algebras, $R_{ij}:A_j\otimes A_i\to A_i\otimes A_j$ twisting maps
for every $i<j$, such that for any $i<j<k$ the maps $R_{ij}$,
$R_{jk}$ and $R_{ik}$ are compatible. Define now for every $i<n-1$
the map
\[ T_{n-1,n}^i:(A_{n-1}\otimes _{R_{n-1\;n}}A_n)\otimes A_i \to
A_i\otimes (A_{n-1}\otimes _{R_{n-1\;n}}A_n)
\]
by $T_{n-1,n}^i:=(R_{i\,n-1}\otimes A_n)\circ (A_{n-1}\otimes
R_{i\,n})$, which are twisting maps for every $i$, as we can
directly apply Theorem \ref{itertwisting} to the maps $R_{i\,n-1}$,
$R_{i\,n}$ and $R_{n-1\,n}$. Furthermore, we have the following
result:

\begin{lemma}
In the above situation, for every $i<j<n-1$, the maps $R_{ij}$,
$T_{n-1,n}^i$ and $T_{n-1,n}^j$ are compatible.
\end{lemma}

\begin{proof}
Using braiding notation the proof can be written as:
\[
\includegraphics[scale=1]{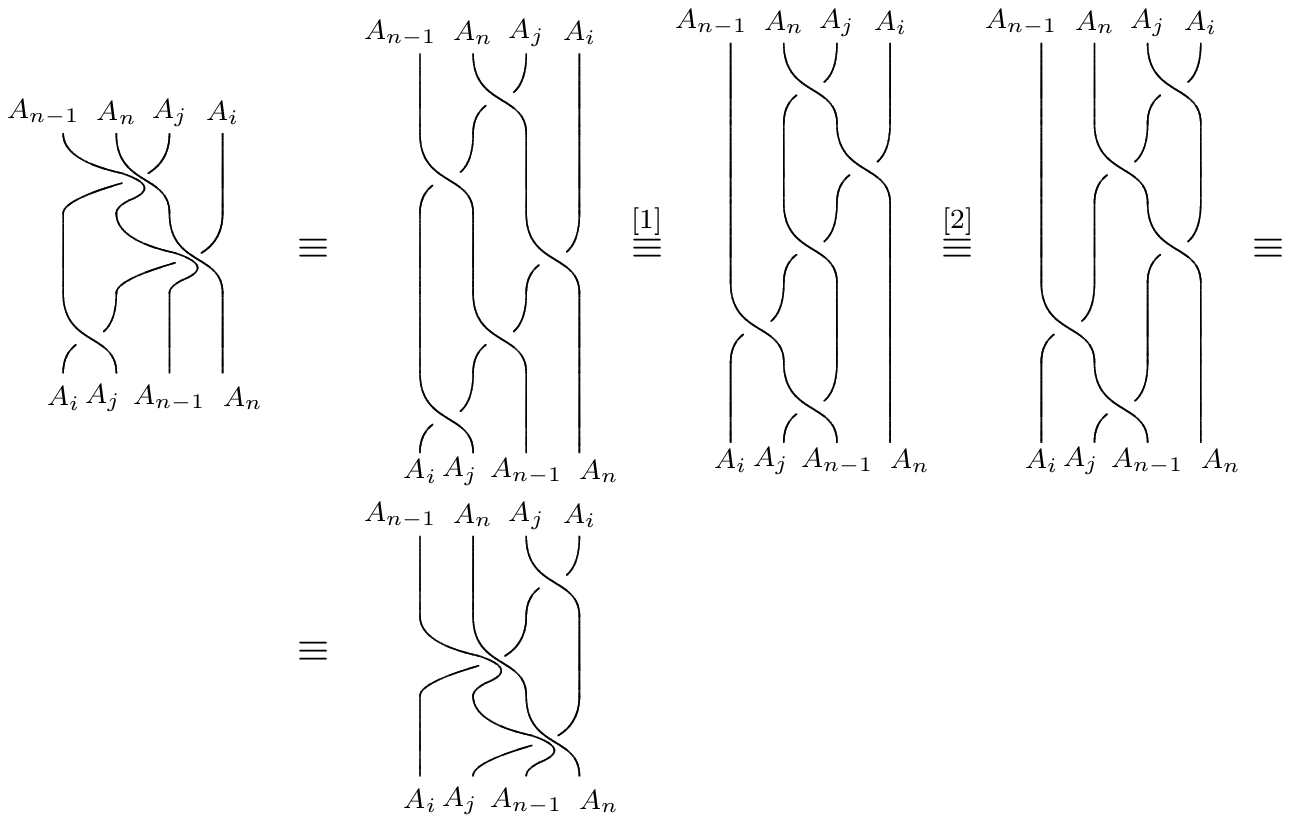}
\]
where in $[1]$ we use the compatibility condition
for $R_{ij}$, $R_{i\,n-1}$ and $R_{j\,n-1}$, and in $[2]$ we use the
compatibility condition for $R_{ij}$, $R_{in}$ and $R_{jn}$.
\qed\end{proof}

So we can apply the induction hypothesis to the $n-1$ algebras
$A_1,\dotsc,A_{n-2}$, and $(A_{n-1}\otimes_{R_{n-1\,n}} A_n)$, and
we obtain that we can build the twisted product of these $n-1$
factors without worrying about parentheses, so we can build the
algebra
\[ A_1\otimes_{R_{12}}\dotsb\otimes
A_{n-2}\otimes_{T_{n-1,n}^{n-2}}(A_{n-1}\otimes_{R_{n-1\,n}}A_n).
\]
Simply observing that
\[A_{n-2}\otimes_{T_{n-1,n}^{n-2}}(A_{n-1}\otimes_{R_{n-1\,n}}A_n)=
(A_{n-2}\otimes_{R_{n-2\,n-1}}A_{n-1})\otimes_{T_{n-2\,n-1}^n}A_n,\]
we see that we could have grouped together any two consecutive
factors. Summarizing, we have sketched the proof of the following
theorem (which we will not write formally to avoid the cumbersome
notation it would involve):

\begin{theorem}[Coherence Theorem]
Let $A_1,\dotsc,A_n$ be algebras, $R_{ij}:A_j\otimes A_i\to
A_i\otimes A_j$ (unital) twisting maps for every $i<j$, such that
for any $i<j<k$ the maps $R_{ij}$, $R_{jk}$ and $R_{ik}$ are
compatible. Then the maps
\[T_{j-1,j}^i:(A_{j-1}\otimes _{R_{j-1\;j}}A_j)\otimes A_i \to A_i\otimes
(A_{j-1}\otimes _{R_{j-1\;j}}A_j)
\]
defined for every $i< j-1$ by $T_{j-1,j}^i:=(R_{i\,j-1}\otimes
A_j)\circ (A_{j-1}\otimes R_{i\,j})$, and the maps
\[T_{j-1,j}^i:A_i\otimes (A_{j-1}\otimes _{R_{j-1\;j}}A_j) \to (A_{j-1}
\otimes _{R_{j-1\;j}}A_j)\otimes A_i
\]
defined for every $i>j$ by $T_{j-1,j}^i:=(A_{j-1}\otimes
R_{j\,i})\circ (R_{j-1\,i}\otimes A_j)$, are twisting maps with the
property that for every $i,k\notin\{j-1,j\}$ the maps $R_{ik}$,
$T_{n-1,n}^i$ and $T_{n-1,n}^k$ are compatible. Moreover, for any
$i$ the (inductively defined) algebras
\[ A_1\otimes_{R_{12}}\dotsb \otimes_{R_{i-3\,i-2}} A_{i-2}\otimes_{T_{i-1,i}^{i-2}}(A_{i-1}\otimes_{R_{i-1\,i}}
A_i)\otimes_{T_{i-1,i}^{i+1}} A_{i+1}\otimes_{R_{i+1\,i+2}}
\dotsb\otimes_{R_{n-1\,n}} A_n
\]
are all equal.
\end{theorem}

As a consequence of this theorem, any property that can be lifted to
iterated twisted tensor products of three factors can be lifted to
products of any number of factors. One of the most interesting
consequences of the Coherence Theorem, or more accurately, of the
former lemma, is that we can state a universal property, analogous
to Theorems \ref{twistuniversalproperty} and \ref{iterunivprop3}. In
order to state the result it is convenient to introduce some
notation. Let us first define the maps
\begin{gather*}
\mathcal{T}_1:A_n\ot \dotsb \ot A_1 \lto A_1\ot A_n\ot \dotsb \ot A_2,\\
\mathcal{T}_1:=(R_{1\,n}\ot \Id_{A_{n-1}\ot \dotsb \ot
A_2})\circ\dotsb\circ (\Id_{A_n\ot \dotsb\ot A_3}\ot R_{12}),\\
\mathcal{T}_2:A_1\ot A_n\ot \dotsb \ot A_2 \lto A_1\ot A_2\ot A_n\ot \dotsb \ot A_3,\\
\mathcal{T}_2:=(A_1\ot R_{2\,n}\ot \Id_{A_{n-1}\ot \dotsb \ot
A_3})\circ\dotsb\circ (\Id_{A_1\ot A_n\ot \dotsb\ot A_4}\ot R_{23}),\\
\vdots\\
\mathcal{T}_{n-1}:A_1\ot \dotsb \ot A_{n-2}\ot A_n\ot A_{n-1} \lto
A_1\ot \dotsb \ot A_{n-2}\ot A_{n-1}\ot A_{n},\\
\mathcal{T}_{n-1}:=A_1\ot \dotsb\ot A_{n-2}\ot R_{n-1\, n},
\end{gather*}
and now define the map
\begin{gather*}
\mathcal{S}: A_n\ot A_{n-1}\ot \dotsb\ot A_1 \lto A_1\ot A_2\ot
\dotsb \ot A_n,\\
\mathcal{S}:=\mathcal{T}_{n-1}\circ\dotsb\circ
\mathcal{T}_2\circ\mathcal{T}_1.
\end{gather*}
With this notation, we can state the Universal Property for iterated
twisted tensor products as follows:

\begin{theorem}[Universal Property]
Let $A_1,\dotsc,A_n$ be algebras, $R_{ij}:A_j\otimes A_i\to
A_i\otimes A_j$ (unital) twisting maps for every $i<j$, such that
for any $i<j<k$ the maps $R_{ij}$, $R_{jk}$ and $R_{ik}$ are
compatible. Suppose that we have an algebra $X$ together with $n$
algebra morphisms $u_i:A_i\to X$ such that
    \begin{equation}
        \mu_X\circ(u_n\ot\dotsb\ot u_1) = \mu_X\circ (u_1\ot \dotsb
        \ot u_n)\circ \mathcal{S}.
    \end{equation}
Then there exists a unique algebra morphism
\[
\vphi:A_1\ot_{R_{12}}A_2 \ot_{R_{23}}\dotsb \ot_{R_{n-1\,n}} A_n
\lto X
\]
such that
\[ \vphi \circ i_{A_j} = u_j,\quad \text{for all $j=1,\dotsc,n$.}
\]
\end{theorem}
\begin{proof}
Following the same procedure as in the proof of Theorem
\ref{iterunivprop3}, it is easy to see that any map $\vphi$
verifying the conditions of the theorem must satisfy
\[ \vphi (a_1\ot\dotsb\ot a_n) = u_1(a_1)\cdot\dotsb \cdot u_n(a_n),
\]
and hence it must be unique. Whence it suffices to define $\vphi$ as
above. The checking of the multiplicative property is also similar
to the one done in the proof of Theorem \ref{iterunivprop3}, and
thus is left to the reader. \qed\end{proof}

A further step in the study of the iterated twisted tensor products
is the lifting of module structures on the factors. Again, if we
have $M$ a left $A$--module, $N$ a left $B$--module, and $P$ a left
$C$--module, the natural way in order to define a left
$(A\otimes_{R_1} B\otimes_{R_2} C)$--module structure on $M\otimes
N\otimes P$ is looking for module twisting maps $\tau_{M,C}:C\otimes
M\to M\otimes C$, $\tau_{M,B}:B\otimes M\to M\otimes B$ and
$\tau_{N,C}:C\otimes N\to N\otimes C$, and defining
    \[ \l_{M\otimes N\otimes P}:=(\l_M\otimes \l_N\otimes \l_P)\circ
    (A\otimes \tau_{M,B}\otimes \tau_{N,C}\otimes P)\circ (A\otimes
    B\otimes \tau_{M,C}\otimes N\otimes P).
    \]
However, as it happened with the iterated product of algebras, in
order to have a left module action it is not enough that
$\tau_{M,C}$, $\tau_{N,C}$ and $\tau_{M,B}$ are module twisting
maps. Realize that, using the $A\otimes_{R_1} B$--module structure
induced on $M\otimes N$ by $\tau_{M,B}$, we can also write the above
action as
    \begin{eqnarray*}
    \l_{M\otimes N\otimes P}&=& (\l_{M\otimes N}\otimes \l_P)\circ
    (A\otimes B\otimes M\otimes \tau_{N,C}\otimes P)\circ (A\otimes
    B\otimes \tau_{M,C} \otimes N\otimes P)\\
    &=& (\l_{M\otimes N}\otimes \l_P)\circ (A\otimes B\otimes \sigma_C\otimes P),
    \end{eqnarray*}
where $\sigma_C: C\otimes (M\otimes N)\to (M\otimes N)\otimes C$ is
defined by $\sigma_C:=(M\otimes \tau_{N,C})\circ(\tau_{M,C}\otimes
N)$, so proving that the three module twisting maps induce a left
module structure on $M\otimes N\otimes P$ is equivalent to prove
that the map $\sigma_C$ is a module twisting map, thus giving a left
$(A\otimes_{R_1}B)\otimes_{T_1} C$--module structure on $(M\otimes
N)\otimes P$. We give sufficient conditions for this to happen in
the following result.

\begin{theorem}
    With the above notation, suppose that the module twisting maps $\tau_{M,C}$, $\tau_{M,B}$
        and  the twisting map $R_2$ satisfy the compatibility
        relation (also called the \dtext{module hexagon condition})
        \begin{equation}\label{modulehexagonequation}
            (M\otimes R_2)\circ (\tau_{M,C}\otimes B)\circ (C\otimes \tau_{M,B})=
            (\tau_{M,B}\otimes C)\circ (B\otimes \tau_{M,C})\circ(R_2\otimes M),
        \end{equation}
        that is, the following diagram
        \[ \xymatrix{
            & C\otimes M\otimes B \ar[r]^{\tau_{M,C}\otimes B}& M\otimes
            C\otimes B \ar[dr]^{M\otimes R_2} \\
            C\otimes B\otimes M \ar[ur]^{C\otimes \tau_{M,B}} \ar[dr]_{R_2\otimes M}
            & & & M\otimes B\otimes C\\
            & B\otimes C\otimes M \ar[r]^{B\otimes \tau_{M,C}}& B\otimes
            M\otimes C \ar[ur]^{\tau_{M,B}\otimes C}
            }
    \]
    is commutative; then:
    \begin{enumerate}
        \item The map $\sigma_C: C\otimes (M\otimes N)\to (M\otimes N)\otimes
        C$ given by $\sigma_C:=(M\otimes \tau_{N,C})\circ(\tau_{M,C}\otimes
N)$ is a module twisting map.
        \item The map $\sigma_{B\otimes C}: (B\otimes C)\otimes M \to M\otimes (B\otimes
        C)$ given by $\sigma_{B\otimes C}:=(\tau_{M,B}\otimes
        C)\circ(B\otimes \tau_{M,C})$ is a module twisting map
        (giving a left $A\otimes_{T_2}(B\otimes_{R_2}C)$--module
        structure on $M\otimes (N\otimes P)$).
    \end{enumerate}
    Moreover, the module structures induced on $M\otimes N\otimes P$ by $\sigma_C$ and $\sigma_{B\otimes
C}$ are equal.
\end{theorem}

\begin{proof}\hfill\\
\noindent\boxed{\textbf{1}}\quad We have to check that $\sigma_C$
satisfies the conditions \eqref{moduletwistcondition1} and
\eqref{moduletwistcondition2}. For the first one, we have that
\[
\includegraphics[scale=1]{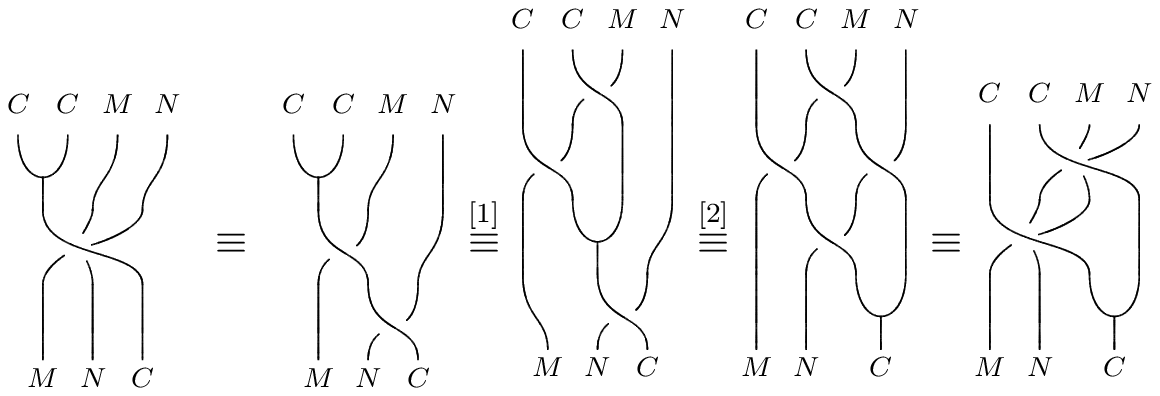}
\]
where in $[1]$ we are using the first module twisting condition for
$\tau_{M,C}$, and in $[2]$ the first module twisting condition for
$\tau_{N,C}$. For the second one, we have
\[
\includegraphics[scale=1]{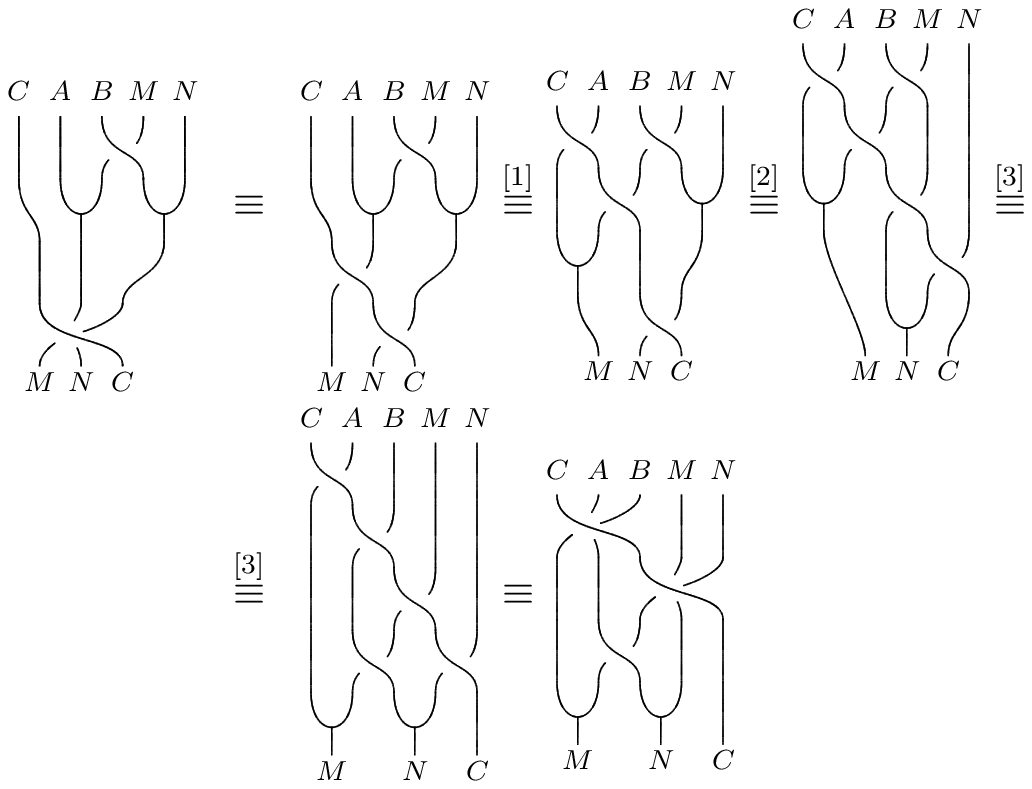}
\]
where in $[1]$ and $[2]$ we use again the module twisting conditions
and in $[3]$ the module hexagon condition.

The proof of $(2)$ is very similar and left to the reader.
\qed\end{proof}

 \begin{remark}
 Note that in this case we cannot prove the module hexagon condition from
 the twisting conditions on the maps. The situation is similar to
 what happens for the case of the existence of module twisting maps.
 It is reasonable to think that some sufficient conditions on the
 modules and the algebras can be given in order to recover the
 converse. For instance, if the modules are free, the situation is
 analogous to the iterated twisting construction for algebras, and
 the converse result can easily be stated.
 \end{remark}


We recall that it is possible to give an explicit description of
modules over various twisted tensor products of algebras arising in
Hopf algebra theory. The same holds in general, as the next result
shows.
\begin{proposition} \label{doua}
Let $A$, $B$ be associative unital algebras and $R:B\ot A\rightarrow
A\ot B$ a unital twisting map. If $M$ is a vector space, then $M$ is
a left $A\ot _RB$-module if and only if it is a left $A$-module and
a left $B$-module satisfying the compatibility condition
    \begin{equation}
        b\cdot (a\cdot m)=a_R\cdot (b_R\cdot m),\ \text{for all $a\in A$, $b\in B$, $m\in
        M$,}
    \end{equation}
where we denote by $\cdot $ both the actions of $A$ and $B$.
\end{proposition}

\begin{proof}
If $M$ is a left $A\ot _RB$-module, it becomes a left $A$-module
with action $a\cdot m=(a\ot 1)\cdot m$ and a left $B$-module with
action $b\cdot m=(1\ot b)\cdot m$. Conversely, it becomes an $A\ot
_RB$-module with action $(a\ot b)\cdot m=a\cdot (b\cdot m)$, details
are left to the reader. \qed\end{proof}

This result can be iterated, generalizing thus the description of
modules over a two-sided smash product from \cite{Panaite02a}.

\begin{proposition} Assume that the hypotheses of Theorem
\ref{itertwisting} are satisfied, such that all algebras and
twisting maps are unital. If $M$ is a vector space, then $M$ is a
left $A\ot _{R_1}B\ot _{R_2}C$-module if and only if it is a left
$A$-module, a left $B$-module, a left $C$-module (all actions are
denoted by $\cdot $) satisfying the compatibility conditions
\begin{eqnarray}
&&b\cdot (a\cdot m)=a_{R_1}\cdot (b_{R_1}\cdot m), \label{compa1} \\
&&c\cdot (b\cdot m)=b_{R_2}\cdot (c_{R_2}\cdot m), \label{compa2} \\
&&c\cdot (a\cdot m)=a_{R_3}\cdot (c_{R_3}\cdot m), \label{compa3}
\end{eqnarray}
for all $a\in A$, $b\in B$, $c\in C$, $m\in M$ (by Proposition
\ref{doua}, these conditions tell that $M$ is a left module over
$A\ot _{R_1}B$, $B\ot _{R_2}C$ and $A\ot _{R_3}C$).
\end{proposition}
\begin{proof}
We only prove that $M$ becomes a left $A\ot _{R_1}B\ot
_{R_2}C$-module with action $(a\ot b\ot c)\cdot m=a\cdot (b\cdot
(c\cdot m))$. We compute (using formula \eqref{trei}):
\begin{eqnarray*}
((a\ot b\ot c)(a'\ot b'\ot c'))\cdot m&=&a(a'_{R_3})_{R_1}\cdot
(b_{R_1}b'_{R_2}\cdot ((c_{R_3})_{R_2}c'\cdot m))\\
\text{\eqref{compa2}}&=&a(a'_{R_3})_{R_1}\cdot (b_{R_1}\cdot
(c_{R_3}\cdot
(b'\cdot (c'\cdot m))))\\
\text{\eqref{compa1}}&=&a\cdot (b\cdot (a'_{R_3}\cdot (c_{R_3}\cdot
(b'\cdot (c'\cdot m)))))\\
\text{\eqref{compa3}}&=&a\cdot (b\cdot (c\cdot (a'\cdot (b'\cdot
(c'\cdot m)))))\\
&=&(a\ot b\ot c)\cdot ((a'\ot b'\ot c')\cdot m),
\end{eqnarray*}
finishing the proof. \qed\end{proof}

Our next result arises as a generalization of the fact from
\cite{Hausser99a}, \cite{BulacuUNa} that a two-sided smash product
over a Hopf algebra is isomorphic to a diagonal crossed product.
\begin{proposition}
Let $(A, B, C, R_1, R_2, R_3)$ be as in Theorem \ref{itertwisting},
and assume that $R_2$ is bijective with inverse $V:B\ot C\rightarrow
C\ot B$. Then $(A, C, B, R_3, V, R_1)$ satisfy also the hypotheses
of Theorem \ref{itertwisting}, and the map $A\ot R_2:A\ot _{R_3}C\ot
_VB\rightarrow A\ot _{R_1}B\ot _{R_2}C$ is an algebra isomorphism.
\end{proposition}
\begin{proof} By Proposition \ref{bij}, $V$ is a twisting map, and
it is obvious that the hexagon condition for $(R_3, V, R_1)$ is
equivalent to the one for $(R_1, R_2, R_3)$. Obviously $A\ot R_2$ is
bijective, we only have to prove that it is an algebra map. This can
be done either by direct computation or, more conceptually, as
follows. Denote $T_2=(R_1\ot C)\circ (B\ot R_3)$ and
$\widetilde{T}_2=(R_3\ot B)\circ (C\ot R_1)$, hence $A\ot _{R_3}C\ot
_VB=A\ot _{\widetilde{T}_2}(C\ot _VB)$ and $A\ot _{R_1}B\ot
_{R_2}C=A\ot _{T_2}(B\ot _{R_2}C)$. By Proposition \ref{bij} we know
that $R_2:C\ot _VB\rightarrow B\ot _{R_2}C$ is an algebra map, and
we obviously have $(A\ot R_2)\circ \widetilde{T}_2=T_2\circ (R_2\ot
A)$, because this is just the hexagon condition. Now it follows that
$A\ot R_2$ is an algebra map, using the following general fact from
\cite{Borowiec00a}: if $A\ot _RB$ and $C\ot _TD$ are twisted tensor
products of algebras and $f:A\rightarrow C$ and $g:B\rightarrow D$
are algebra maps satisfying the condition $(f\ot g)\circ R=T\circ
(g\ot f)$, then $f\ot g:A\ot _RB\rightarrow C\ot _TD$ is an algebra
map. \qed\end{proof}


As our main motivations aimed at applications of our construction to
the field of noncommutative geometry, we are especially interested
in finding processes that allow us to lift constructions associated
to geometrical invariants of the algebras to their (iterated)
twisted tensor products. Among these geometrical invariants, the
first one to be taken into account is of course the algebra of
differential forms. For the case of the twisted product of two
algebras, a twisted product of the algebras of universal
differential forms is build in a unique way, as it is shown in
Theorem \ref{twistdifferentialforms}; there, the construction of
these extended twisting maps is deduced from the universal property
of the first order universal differential calculus. This extension
is compatible with our extra condition for constructing iterated
products, as we show in the following result:

\begin{theorem}\label{itertwistdifferentialforms}
Let $A$, $B$, $C$ be algebras, and let $R_1:B\otimes A\lto A\otimes
B$, $R_2:C\otimes B\lto B\otimes C$, $R_3:C\otimes A\lto A\otimes C$
be twisting maps satisfying the hexagon equation, then the extended
twisting maps $\widetilde{R}_1$, $\widetilde{R}_2$ and
$\widetilde{R}_3$ also satisfy the hexagon equation. Moreover, $\O
A\otimes_{\widetilde{R}_1}\O B\otimes_{\widetilde{R}_2}\O C$ is a
differential graded algebra, with differential
    \[d = d_A\otimes \O B\otimes \O C + \eps_A\otimes d_B\otimes \O C + \eps_A\otimes \eps_B\otimes d_C.
    \]
\end{theorem}

\begin{proof}
For proving that the extended twisting maps satisfy the hexagon
equation, we use a standard technique when dealing with algebras of
differential forms.

Firstly, observe that when restricted to the zero degree part of the
algebras of differential forms, the extended twisting maps coincide
with the original ones, and hence they trivially satisfy the hexagon
equation.

Now, suppose that we have elements $\omega\in \O A$, $\eta\in \O B$,
$\theta\in \O C$ such that the hexagon equation is satisfied when
evaluated on $\omega\otimes \eta\otimes\theta$, and let us show that
then the hexagon equation is also satisfied when evaluated in
$d_A\omega\otimes \eta\otimes\theta$, $\omega\otimes
d_B\eta\otimes\theta$ and $\omega\otimes \eta\otimes d_C\theta$,
that is, we will show that the hexagon condition is stable under
application of any of the differentials $d_A$, $d_B$ and $d_C$.

Let us start proving that the condition holds for $\omega\otimes
\eta\otimes d_C\theta$. Using again braiding notation, we have
\[
\includegraphics[scale=1]{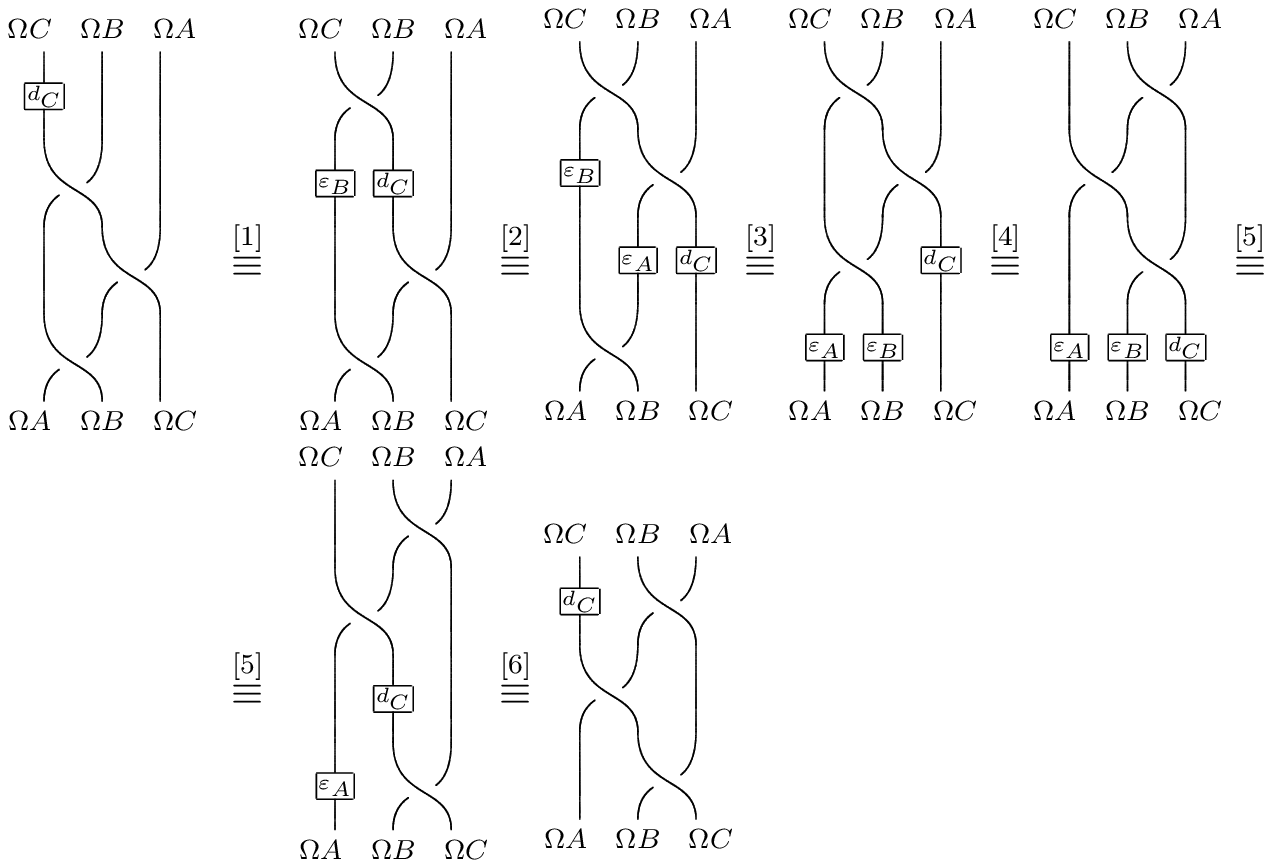}
\]
where in $[1]$, $[2]$, $[5]$ and $[6]$ we are using the property
\eqref{twisteddiff1} for $d_C$ with respect to $R_{2}$ and $R_3$
respectively, in $[3]$ the (obvious) fact that the gradings commute
with the extended twisting maps (since they are homogeneous), and in
$[4]$ we are using the hexagon equation for
$\omega\otimes\eta\otimes \theta$. The corresponding proofs for
$\omega\otimes d\eta\otimes \theta$ and $d\omega\otimes\eta\otimes
\theta$ are almost identical. Summarizing, the hexagon condition is
stable under differentials in $\O A$, $\O B$ and $\O C$.

 Finally, suppose that we have elements $\omega\in \O A$, $\eta\in \O B$, $\theta_1, \theta_2\in \O C$ such that the
hexagon equation is satisfied when evaluated on $\omega\otimes
\eta\otimes\theta_1$ and $\omega_2\otimes \eta\otimes\theta_2$, and
let us show that in this case the hexagon condition also holds on
$\omega\otimes \eta\otimes\theta_1\theta_2$:
\[
\includegraphics[scale=1]{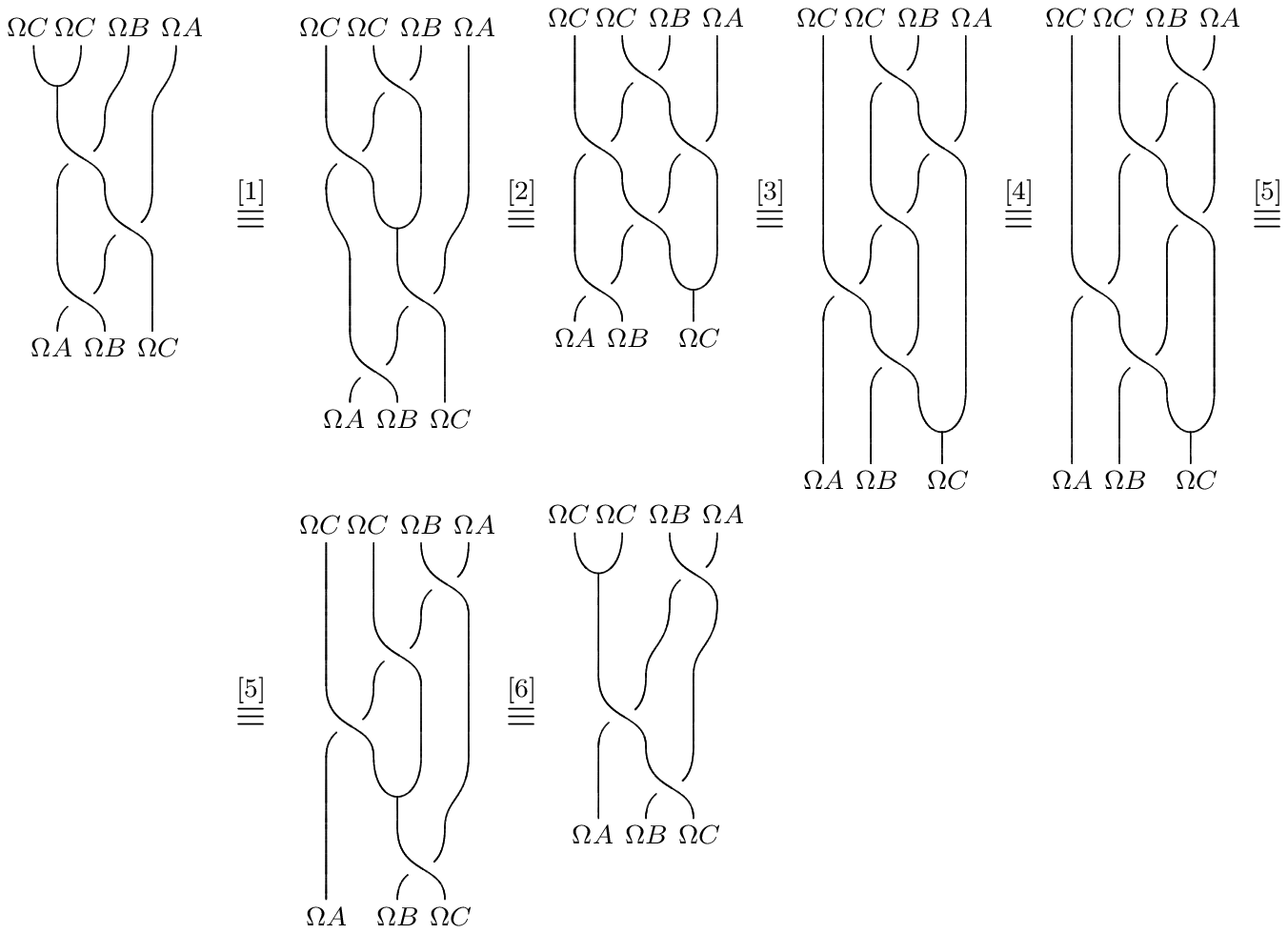}
\]
where in $[1]$, $[2]$, $[5]$ and $[6]$ we use the pentagon equations
\eqref{pentagonequation2} for the twisting maps $\widetilde{R}_2$
and $\widetilde{R}_3$, and in $[3]$ and $[4]$ we use the hexagon
condition for $\omega\otimes\eta\otimes \theta_1$ and
$\omega\otimes\eta\otimes \theta_2$ respectively. In a completely
analogous way we can prove that the hexagon condition holds for
$\omega\otimes\eta_1\eta_2\otimes \theta$ and
$\omega_1\omega_2\otimes\eta\otimes \theta$, that is: the hexagon
condition remains stable under products in $\O A$, $\O B$ and $\O
C$.

Now, taking into account that $\O A$, $\O B$ and $\O C$ are
generated as graded differential algebras by the elements of degree
0, we may conclude that the hexagon condition holds completely.

In order to prove that $\O A\otimes_{\widetilde{R}_1} \O
B\otimes_{\widetilde{R}_2}\O C$ is a graded differential algebra, it
is enough to observe that $\O A\otimes_{\widetilde{R}_1} \O
B\otimes_{\widetilde{R}_2}\O C= (\O A\otimes_{\widetilde{R}_1} \O
B)\otimes_{\widetilde{T}_2}\O C$, the last being (because of Theorem
\ref{twistdifferentialforms}) a  graded differential algebra with
differential
\[ d= d_{A\otimes_{R_1} B}\otimes \O C + \eps_{A\otimes_{R_1}
B}\otimes d_C,
\]
and taking into account that
\begin{gather*}
d_{A\otimes_{R_1} B} = d_A\otimes \O B + \eps_A\otimes
d_B,\\
\eps_{A\otimes_{R_1} B}= \eps_A\otimes \eps_B,
\end{gather*}
we obtain
    \[  d = d_A\otimes \O B\otimes \O C + \eps_A\otimes d_B\otimes \O C + \eps_A\otimes \eps_B\otimes
    d_C,
    \]
as we wanted to show. \qed\end{proof}

As most of our motivation comes from some algebras used in Connes'
theory, in order to deal properly with $\ast$--algebras we would
like to find a suitable extension of condition \eqref{starcondition}
to our framework. As the definition of the involution in a twisted
tensor product also involves the usual flip $\tau$, before extending
the conditions to an iterated product, we need a technical (and easy
to prove) result:

\begin{lemma}\label{twoflipscondition}
Let $A$, $B$, $C$ be algebras, and let $R:B\otimes A\to A\otimes B$
be a twisting map. Consider also the usual flips $\tau_{BC}:B\otimes
C\to C\otimes B$ and $\tau_{AC}:A\otimes C\to C\otimes A$, then the
maps $\tau_{AC}$, $R$ and $\tau_{BC}$ satisfy the hexagon condition
(in $B\otimes A\otimes C$).
\end{lemma}

\begin{proof}
Just write down both sides of the equation and realize they are
equal. \qed\end{proof}

\begin{remark}
In general, we can say that any twisting map is compatible with a
pair of usual flips, regardless the ordering of the factors. As the
inverse of a usual flip is also a usual flip, we may also use this
result when one of the flips is inverted.
\end{remark}

Similarly to what happened with differential forms, in order to be
able to extend the involutions to the iterated product, it is enough
that condition \eqref{starcondition} is satisfied for every pair of
algebras. More concretely, we have the following result:

\begin{theorem}
Let $A$, $B$, $C$ be $\ast$--algebras with involutions $j_A$, $j_B$
and $j_C$ respectively, $R_1:B\otimes A\to A\otimes B$,
$R_2:C\otimes B\to B\otimes C$ and $R_3:C\otimes A\to A\otimes C$
compatible twisting maps such that
\begin{eqnarray}
    (R_1\circ (j_B\otimes j_A)\circ \tau_{AB})\circ (R_1\circ (j_B\otimes j_A)\circ
    \tau_{AB}) & = & A\otimes B,\label{involution1} \\
    (R_2\circ (j_C\otimes j_B)\circ \tau_{BC})\circ (R_2\circ (j_C\otimes
    j_B)\circ \tau_{BC}) & = & B\otimes C, \label{involution2} \\
    (R_3\circ (j_C\otimes j_A)\circ \tau_{AC})\circ (R_3\circ (j_C\otimes j_A)\circ
    \tau_{AC}) & = & A\otimes C.\label{involution3}
\end{eqnarray}
Then $A\otimes_{R_1} B\otimes_{R_2} C$ is a $\ast$--algebra with
involution
\[ j=(R_1\otimes C)\circ (B\otimes R_3)\circ (R_2\otimes
A)\circ(j_C\otimes j_B\otimes j_A)\circ (C\otimes \tau_{AB})\circ
(\tau_{AC}\otimes B)\circ (A\otimes\tau_{BC}),
\]
where $\tau_{AB}:A\otimes B\to B\otimes A$, $\tau_{BC}:B\otimes C\to
C\otimes B$, and $\tau_{AC}:A\otimes C\to C\otimes A$ denote the
usual flips.
\end{theorem}

\begin{proof}
Consider $j$ defined as above, and let us check that it is an
involution, i. e., that $j^2=A\otimes B\otimes C$. Firstly, observe
that, if we denote by $\tau$ all the usual flips and by $\bar{\tau}$
their inverses, we have that
\[
\includegraphics[scale=1]{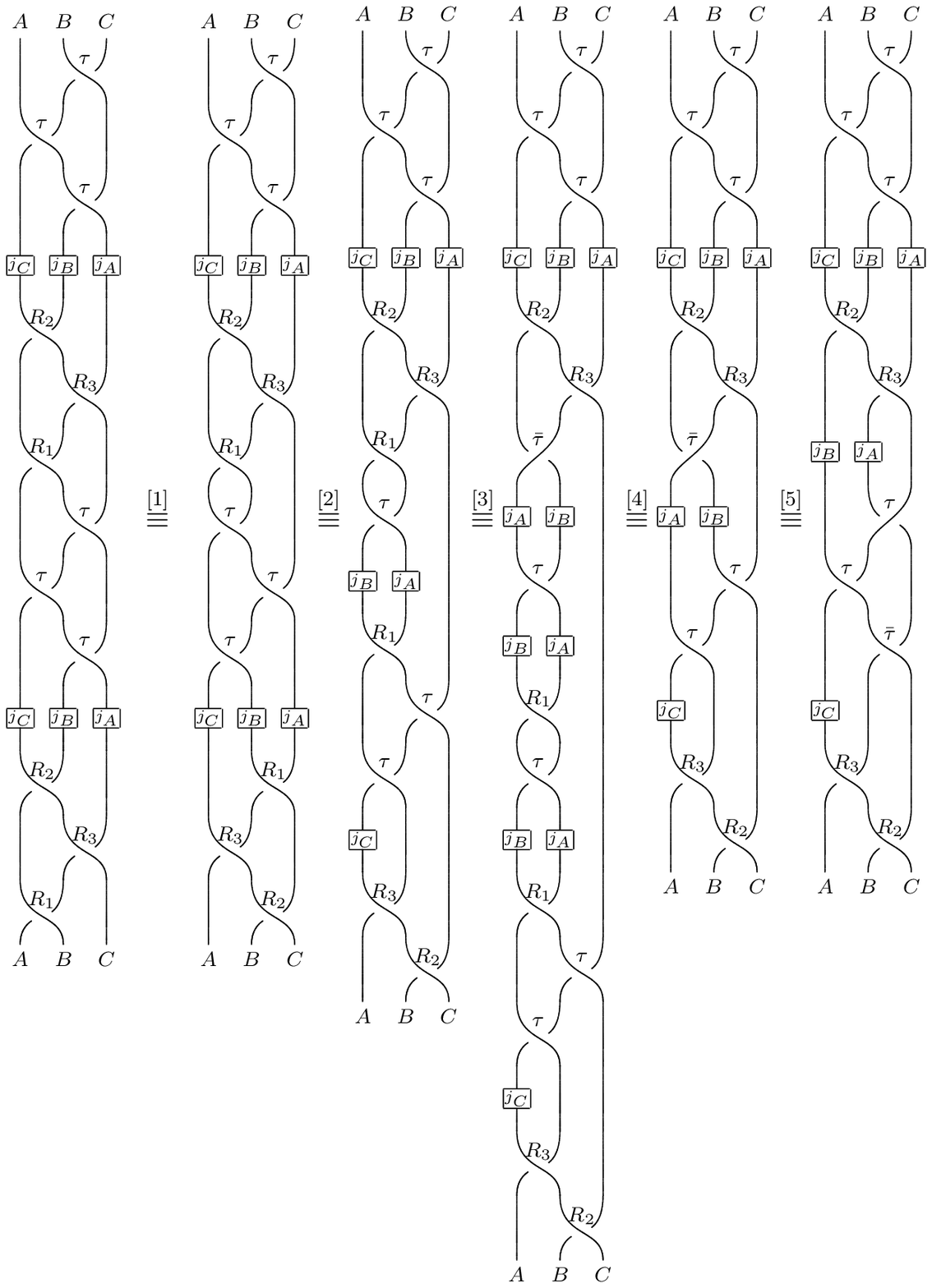}
\]
where in $[1]$ we use the hexagon conditions for the flips (which is
obvious) and the hexagon conditions for $R_1$, $R_2$, $R_3$, in
$[2]$ we use the fact that the involutions $j_A$ and $j_B$ commute
with the flips, and the hexagon condition for $R_1$ and two flips
(as stated in the former lemma). Equivalence $[3]$ is due to the
fact that both the square of the involutions, and the composition of
a flip with its inverse are the identity. In $[4]$ we apply
\eqref{involution1}, and in $[5]$ we use again that the involutions
commute with the flips, plus the hexagon condition for
$\tau_{AB}^{-1}$ and two usual flips. To conclude the proof, observe
that \clearpage
\[
\includegraphics[scale=1]{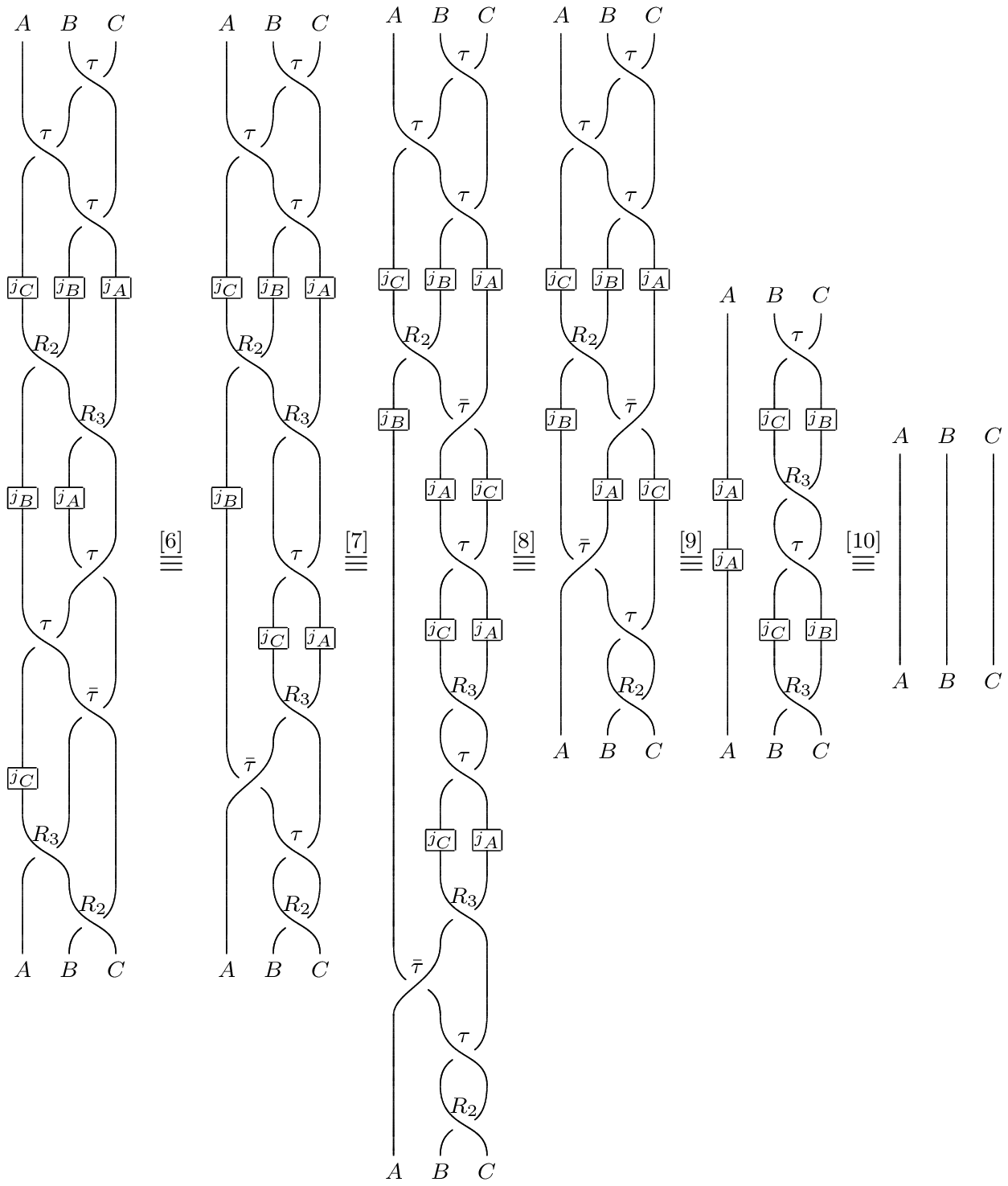}
\]
where in $[6]$ we apply (twice) the commutation of $j_C$ with the
flips, plus the hexagon for $R_3$ and two flips (again because of
the former lemma). Equality $[7]$ holds again because we are just
adding a term (two squared involutions, a flip, and its inverse)
that equals the identity, while $[8]$ holds by applying
\eqref{involution2}. $[9]$ is due to the fact that in the last
diagram the element of $A$ is not modified at all, since all the
crossings are usual flips, and we get $[10]$ using
\eqref{involution2} and the fact that $j_A$ is an involution.
\qed\end{proof}


\section{Examples}\label{sec:examples}
\subsection{Generalized smash products}
We begin by recalling the construction of the so-called generalized
smash products. Let $H$ be a bialgebra. For a right $H$-comodule
algebra $({\mf A}, \r )$ we denote $\r (\mfa )=\mfa _{<0>}\ot \mfa
_{<1>}$, for any $\mfa \in {\mf A}$. Similarly, for a left
$H$-comodule algebra $({\mf B}, \l )$, if $\mfb \in {\mf B}$ then we
denote $\l (\mfb )=\mfb _{[-1]}\ot \mfb _{[0]}$.

Let $A$ be a left $H$-module algebra and $\mathfrak{B}$ a left
$H$-comodule algebra. Denote by $A\gsm \mathfrak{B}$ the $k$-vector
space $A\ot\mathfrak{B}$ with newly defined multiplication
\begin{equation}\label{gsm}
(a\gsm \mf {b})(a'\gsm \mf {b}')= a(\mf {b}_{[-1]}\cd a')\gsm \mf
{b}_{[0]}\mf {b}' ,
\end{equation}
for all $a, a'\in A$ and $\mf {b}, \mf {b}'\in \mathfrak{B}$. Then
$A\gsm \mathfrak{B}$ is an associative algebra with unit $1_A\gsm
1_{\mathfrak{B}}$. If we take $\mathfrak{B}=H$ then $A\gsm H$ is
just the ordinary smash product $A\# H$, whose multiplication is
\begin{eqnarray*}
&&(a\# h)(a'\# h')=a(h_1\cdot a')\# h_2h'.
\end{eqnarray*}
The algebra $A\gsm \mathfrak{B}$ is called the (left) generalized
smash product
of $A$ and $\mathfrak{B}$. \\
Similarly, if $B$ is a right $H$-module algebra and $\mathfrak{A}$
is a right $H$-comodule algebra, then we denote by $\mf {A}\gtl B$
the $k$-vector space $\mf {A}\ot B$ with the newly defined
multiplication
\begin{equation}\label{gtl}
(\mf {a}\gtl b)(\mf {a}{'}\gtl b{'})= \mf {a}\mf {a}'_{<0>}\gtl
(b\cd \mf {a}'_{<1>})b{'},
\end{equation}
for all $\mf {a}, \mf {a}{'}\in {\mf A}$ and $b, b{'}\in B$. Then
$\mf {A}\gtl B$ is an associative algebra with unit $1_{\mf {A}}\gtl
1_B$, called also the (right) generalized smash product of $\mf {A}$
and $B$.

We recall some facts from \cite{BulacuUNa}. Let $H$ be a bialgebra,
$A$ a left $H$-module algebra, $B$ a right $H$-module algebra and
${\mb A}$ an $H$-bicomodule algebra. Then $A\gsm {\mb A}$ becomes a
right $H$-comodule algebra with structure
\begin{eqnarray*}
&&A\gsm {\mb A}\rightarrow (A\gsm {\mb A})\ot H, \;\;\; a\gsm
u\mapsto (a\gsm u_{<0>})\ot u_{<1>},
\end{eqnarray*}
and ${\mb A}\gtl B$ becomes a left $H$-comodule algebra with
structure
\begin{eqnarray*}
&&{\mb A}\gtl B\rightarrow H\ot ({\mb A}\gtl B), \;\;\; u\gtl
b\mapsto u_{[-1]}\ot (u_{[0]}\gtl b).
\end{eqnarray*}
Moreover, we have:
\begin{proposition} (\cite{BulacuUNa}) $(A\gsm {\mb A})\gtl B\equiv
A\gsm ({\mb A}\gtl B)$ as algebras. If ${\mb A}=H$, this algebra is
denoted by $A\# H\# B$ and is called a two-sided smash product.
\end{proposition}
This result is a particular case of Theorem \ref{itertwisting}.
Indeed, define the maps
\begin{eqnarray*}
&&R_1:{\mb A}\ot A\rightarrow A\ot {\mb A}, \;\;\;R_1(u\ot a)=
u_{[-1]}\cdot a\ot u_{[0]}, \\
&&R_2:B\ot {\mb A}\rightarrow {\mb A}\ot B, \;\;\;R_2(b\ot
u)=u_{<0>}\ot
b\cdot u_{<1>}, \\
&&R_3:B\ot A\rightarrow A\ot B, \;\;\;R_3(b\ot a)=a\ot b,
\end{eqnarray*}
which obviously are twisting maps because $A\ot _{R_1}{\mb A}=A\gsm
{\mb A}$ and ${\mb A}\ot _{R_2}B={\mb A}\gtl B$ are associative
algebras. Moreover, if we define the maps
\begin{eqnarray*}
&&T_1:B\ot (A\ot {\mb A})\rightarrow (A\ot {\mb A})\ot B, \;\;\;
T_1:=(A\ot R_2)\circ (R_3\ot {\mb A}), \\
&&T_2:({\mb A}\ot B)\ot A\rightarrow A\ot ({\mb A}\ot B), \;\;\;
T_2:=(R_1\ot B)\circ ({\mb A}\ot R_3),
\end{eqnarray*}
then one can see that
\begin{eqnarray*}
&&(A\gsm {\mb A})\ot _{T_1}B=(A\gsm {\mb A}) \gtl B, \;\;\;A\ot
_{T_2}({\mb A}\gtl B)=A\gsm ({\mb A}\gtl B).
\end{eqnarray*}
\subsection{Generalized diagonal crossed products}
We recall the construction of the so-called generalized diagonal
crossed product, cf. \cite{BulacuUNa}, \cite{Hausser99a}. Let $H$ be
a Hopf algebra with bijective antipode $S$, ${\mathcal{A}}$ an
$H$-bimodule algebra and $\mb {A}$ an $H$-bicomodule algebra. Then
the generalized diagonal crossed product ${\mathcal{A}}\bowtie \mb
{A}$ is the following associative algebra structure on
${\mathcal{A}}\ot {\mb A}$:
\begin{equation}\label{gdphopf}
(\varphi \bowtie u)(\varphi {'}\bowtie u{'})= \varphi (u_{\{-1\}}\cd
\varphi {'}\cd \smi (u_{\{1\}})) \bowtie u_{\{0\}}u{'},
\end{equation}
for all $\varphi , \varphi {'}\in {\mathcal{A}}$ and $u, u{'}\in \mb
{A}$, where
\begin{eqnarray*}
&&u_{\{-1\}}\ot u_{\{0\}}\ot u_{\{1\}} :=u_{<0>_{[-1]}}\ot
u_{<0>_{[0]}}\ot u_{<1>}=u_{[-1]}\ot u_{[0]_{<0>}}\ot u_{[0]_{<1>}}.
\end{eqnarray*}

We recall some facts from \cite{PanaiteUNa}. Let $H$ be a Hopf
algebra with bijective antipode $S$, ${\mathcal{A}}$ an $H$-bimodule
algebra and ${\mb A}$ an $H$-bicomodule algebra. Let also $A$ be an
algebra in the Yetter-Drinfeld category $_H^H{\mathcal{YD}}$, that
is, $A$ is a left $H$-module algebra, a left $H$-comodule algebra
(with left $H$-comodule structure denoted by $a\mapsto a_{(-1)}\ot
a_{(0)}\in H\ot A$) and the Yetter-Drinfeld compatibility condition
holds:
\begin{eqnarray}
&&h_1a_{(-1)}\ot h_2\cdot a_{(0)}=(h_1\cdot a)_{(-1)}h_2\ot
(h_1\cdot a)_{(0)}, \;\;\;\forall \;h\in H, \;a\in A. \label{yd}
\end{eqnarray}
Consider first the generalized smash product ${\mathcal{A}}\gsm A$,
as associative algebra. From the condition (\ref{yd}), it follows
that ${\mathcal{A}}\gsm A$ becomes an $H$-bimodule algebra, with
$H$-actions
\begin{eqnarray*}
&&h\cdot (\varphi \gsm a)=h_1\cdot \varphi \gsm h_2\cdot a, \\
&&(\varphi \gsm a)\cdot h=\varphi \cdot h\gsm a,
\end{eqnarray*}
for all $h\in H$, $\varphi \in {\mathcal{A}}$ and $a\in A$, hence we
may consider
the algebra $({\mathcal{A}}\gsm A)\bowtie {\mb A}$.\\
Then, consider the generalized smash product $A\gsm {\mb A}$, as
associative algebra. Using the condition (\ref{yd}), one can see
that $A\gsm {\mb A}$ becomes an $H$-bicomodule algebra, with
$H$-coactions
\begin{eqnarray*}
&&\rho :A\gsm {\mb A}\rightarrow (A\gsm {\mb A})\ot H,\;\;\;\rho
(a\gsm u)=
(a\gsm u_{<0>})\ot u_{<1>}, \\
&&\lambda :A\gsm {\mb A}\rightarrow H\ot (A\gsm {\mb A}), \;\;\;
\lambda (a\gsm u)=a_{(-1)}u_{[-1]}\ot (a_{(0)}\gsm u_{[0]}),
\end{eqnarray*}
for all $a\in A$ and $u\in {\mb A}$, hence we may consider the
algebra ${\mathcal{A}}\bowtie (A\gsm {\mb A})$.

A similar computation to the one in the proof of Proposition 3.4 in
\cite{PanaiteUNa} shows:
\begin{proposition}\label{cucu}
We have an algebra isomorphism $({\mathcal{A}}\gsm A)\bowtie {\mb
A}\equiv {\mathcal{A}}\bowtie (A\gsm {\mb A})$, given by the trivial
identification.
\end{proposition}
This result is also a particular case of Theorem \ref{itertwisting}.
Indeed, define the maps:
\begin{eqnarray*}
&&R_1:A\ot {\mathcal{A}}\rightarrow {\mathcal{A}}\ot A,
\;\;\;R_1(a\ot \varphi )=
a_{(-1)}\cdot \varphi \ot a_{(0)},\\
&&R_2:{\mb A}\ot A\rightarrow A\ot {\mb A}, \;\;\;R_2(u\ot a)=
u_{[-1]}\cdot a\ot u_{[0]}, \\
&&R_3:{\mb A}\ot {\mathcal{A}}\rightarrow {\mathcal{A}}\ot {\mb
A},\;\;\; R_3(u\ot \varphi )=u_{\{-1\}}\cdot \varphi \cdot
S^{-1}(u_{\{1\}})\ot u_{\{0\}},
\end{eqnarray*}
which are all twisting maps because ${\mathcal{A}}\ot
_{R_1}A={\mathcal{A}}\gsm A$, $A\ot _{R_2}{\mb A}=A\gsm {\mb A}$ and
${\mathcal{A}}\ot _{R_3}{\mb A}= {\mathcal{A}}\bowtie {\mb A}$ are
associative algebras. Now, if we define the maps
\begin{eqnarray*}
&&T_1:{\mb A}\ot ({\mathcal{A}}\ot A)\rightarrow ({\mathcal{A}}\ot
A)\ot {\mb A}, \;\;\;
T_1:=({\mathcal{A}}\ot R_2)\circ (R_3\ot A), \\
&&T_2:(A\ot {\mb A})\ot {\mathcal{A}}\rightarrow {\mathcal{A}}\ot
(A\ot {\mb A}), \;\;\; T_2:=(R_1\ot {\mb A})\circ (A\ot R_3),
\end{eqnarray*}
then one can check that we have
\begin{eqnarray*}
&&({\mathcal{A}}\gsm A)\ot _{T_1}{\mb A}=({\mathcal{A}}\gsm
A)\bowtie {\mb A}, \;\;\; {\mathcal{A}}\ot _{T_2}(A\gsm {\mb
A})={\mathcal{A}}\bowtie (A\gsm {\mb A}),
\end{eqnarray*}
hence indeed we recover Proposition \ref{cucu}.

\subsection{The noncommutative $2n$--planes}
Recall from section \ref{sec:prelim} that the noncommutative plane
associated to an antisymmetric matrix, $\theta=(\theta_{\mu\nu}) \in
\M_n(\R)$, is the associative algebra $C_{alg}(\R^{2n}_\theta)$
generated by $2n$ elements $\{z^\mu,\bar{z}^\mu\}_{\mu=1,\dotsc,n}$
with relations
\begin{equation*}
\left.\begin{array}{r}
z^\mu z^\nu=\l^{\mu\nu}z^\nu z^\mu \\
\bar{z}^\mu \bar{z}^\nu=\l^{\mu\nu}\bar{z}^\nu \bar{z}^\mu \\
\bar{z}^\mu z^\nu = \l^{\nu \mu} z^\nu \bar{z}^\mu
\end{array}\right\}\forall\,\mu,\nu=1,\dotsc,n,\ \text{being
$\l^{\mu\nu}:=e^{i\theta_{\mu\nu}}$,}
\end{equation*}
and endowed with the $\ast$--operation induced by
$(z^\mu)^\ast:=\bar{z}^\mu$ (cf. \cite{Connes02a}).

Observe that as $\theta$ is antisymmetric, we have that
$z^\mu\bar{z}^\mu=\bar{z}^\mu z^\mu$, so for every $\mu=1,\dotsc,n$
the algebra $A_\mu$ generated by the elements $z^\mu$ and
$\bar{z}^\mu$ is commutative, so $A_\mu\cong \c[z^\mu,\bar{z}^\mu]$.
We have then $n$ commutative algebras (indeed, $n$ copies of the
polynomial algebra in two variables) contained in the noncommutative
plane. Consider, for $\mu<\nu$, the mappings defined on generators
by
\begin{eqnarray*}
R_{\mu\nu}:\c[z^\nu,\bar{z}^\nu]\otimes \c[z^\mu,\bar{z}^\mu] & \lto
&
\c[z^\mu,\bar{z}^\mu]\otimes \c[z^\nu,\bar{z}^\nu],\\
z^\nu\otimes z^\mu & \lmto & \l^{\nu\mu}z^\mu\otimes z^\nu,\\
\bar{z}^\nu\otimes \bar{z}^\mu & \lmto & \l^{\nu\mu}\bar{z}^\mu\otimes \bar{z}^\nu,\\
\bar{z}^\nu\otimes z^\mu & \lmto & \l^{\mu\nu}z^\mu\otimes \bar{z}^\nu,\\
z^\nu\otimes \bar{z}^\mu & \lmto & \l^{\mu\nu}\bar{z}^\mu\otimes
z^\nu.
\end{eqnarray*}
Obviously these formulae extend in a unique way to (unital) twisting
maps $R_{\mu\nu}$. Condition \eqref{starcondition} is trivially
satisfied, so every possible twisted tensor product is still a
$\ast$--algebra. As on the algebra generators our twisting map is
just the usual flip multiplied by a constant, the hexagon condition
is also satisfied in a straightforward way. The iterated twisted
tensor product
\[ C[z^1,\bar{z}^1]\otimes_{R_{12}} C[z^2,\bar{z}^2] \otimes_{R_{23}} \dotsb
\otimes_{R_{n-1\, n}} C[z^n,\bar{z}^n]
\]
is isomorphic to the noncommutative $2n$--plane
$C_{alg}(\R^{2n}_\theta)$. Furthermore, for every $\mu=1,\dotsc,n$,
let $\O_\mu := \O_{alg}(\R^2)$ be the graded differential algebra of
algebraic differential forms build over the algebra
$\c[z^\mu,\bar{z}^\mu]$, and observe that for $\mu < \nu$ the map
$\overline{R}_{\mu\nu}:\O_\nu\otimes \O_\mu \lto \O_\mu\otimes
\O_\nu$ defined on generators by
\[
\begin{array}{rclcrcl}
    z^\nu\otimes z^\mu & \lmto & \l^{\nu\mu}z^\mu\otimes z^\nu,
    &\qquad&
    \bar{z}^\nu\otimes \bar{z}^\mu & \lmto & \l^{\nu\mu}\bar{z}^\mu\otimes \bar{z}^\nu,\\
    \bar{z}^\nu\otimes z^\mu & \lmto & \l^{\mu\nu}z^\mu\otimes
    \bar{z}^\nu, &&
    z^\nu\otimes \bar{z}^\mu & \lmto & \l^{\mu\nu}\bar{z}^\mu\otimes z^\nu,\\
    dz^\nu \otimes dz^\mu & \lmto & -\l^{\nu\mu}dz^\mu \otimes
    dz^\nu,&&
    d\bar{z}^\nu \ot d\bar{z}^\mu & \lmto & -\l^{\nu\mu}d\bar{z}^\mu \ot d\bar{z}^\nu,\\
    d\bar{z}^\nu \ot dz^\mu & \lmto &   -\l^{\mu\nu}dz^\mu \ot
    d\bar{z}^\nu,&&
    dz^\nu \ot d\bar{z}^\mu & \lmto &  -\l^{\mu\nu}d\bar{z}^\mu \ot dz^\nu,\\
    z^\nu \ot dz^\mu & \lmto & \l^{\nu\mu}dz^\mu \ot z^\nu,&&
    \bar{z}^\nu \ot d\bar{z}^\mu & \lmto & \l^{\nu\mu}d\bar{z}^\mu \ot \bar{z}^\nu,\\
    \bar{z}^\nu \ot dz^\mu & \lmto & \l^{\mu\nu}dz^\mu \ot
    \bar{z}^\nu,&&
    z^\nu \ot d\bar{z}^\mu & \lmto & \l^{\mu\nu}d\bar{z}^\mu \ot z^\nu,
\end{array}
\]
extends in a unique way to a twisting map defined on $\O_\nu\otimes
\O_\mu$. This twisting map satisfies conditions \eqref{twisteddiff1}
and \eqref{twisteddiff2}, hence, by the uniqueness of the twisting
map extension to the algebras of differential forms given by Theorem
\ref{twistdifferentialforms}, the maps $\overline{R}_{\mu\nu}$
coincide with the maps $\widetilde{R}_{\mu\nu}$ obtained in the
theorem. So, by applying Theorem \ref{itertwistdifferentialforms} it
follows that they are compatible twisting maps. It is then easy to
check that the iterated twisted tensor product
$\O_1\ot_{\overline{R}_{12}} \dotsb\ot_{\overline{R}_{n-1\,n}}\O_n$
is isomorphic, as a graded (involutive) differential algebra, to the
algebra $\O_{alg}(\R^{2n}_\theta)$ of algebraic differential forms
on the noncommutative $2n$--plane.

\subsection{The Observable Algebra of Nill--Szlachányi}
In \cite{Nill97a}, Nill and Szlachányi construct, given a finite
dimensional $C^\ast$--Hopf algebra $H$ and its dual $\hat{H}$, the
\dtext{algebra of observables}, denoted by $\mathcal{A}$, by means
of the smash products defined by the natural actions existing
between $H$ and $\hat{H}$. Their interest in studying such an
algebra arises as it turns out to be the \emph{observable algebra of
a generalized quantum spin chain with $H$--order and
$\hat{H}$--disorder symmetries}, and they also observe that when
$H=\c G$ is a group algebra this algebra $\mathcal{A}$ becomes an
ordinary $G$--spin model. We do not need here the physical
interpretation of this algebra, our aim is to show that the
construction of the algebra $\mathcal{A}$ carried out in
\cite{Nill97a} fits inside our framework of iterated twisted tensor
products.

We start by fixing $H$ a finite dimensional $C^\ast$--Hopf algebra,
that is, a $C^\ast$--algebra endowed with a comultiplication
$\Delta:H\to H\otimes H$, a counit $\eps:H\to \c$ and an antipode
$S:H\to H$ satisfying the usual compatibility relations required for
defining Hopf algebras, and with the extra assumptions that $\Delta$
and $\eps$ are $\ast$--algebra morphisms, and such that
$S(S(x)^\ast)^\ast=x$ for all $x\in H$ (see \cite[Section
IV.8]{Kassel95a} for details). If $H$ is a $\ast$--Hopf algebra, it
follows that $S^{-1}=\bar{S}=\ast\circ S\circ \ast$ is the antipode
of the opposite Hopf algebra $H_{\text{op}}$ (see \cite{Sweedler69a}
for details). The dual Hopf algebra of a $\ast$--Hopf algebra is
also a $\ast$--Hopf algebra, with involution given by
$\vphi^{\ast}:=S(\vphi_\ast)$, where $\vphi\mapsto \vphi_\ast$ is
the antilinear involutive algebra automorphism given by
$\vphi_\ast(x):=\overline{\vphi(x^\ast)}$. We have canonical
pairings between $H$ and $\hat{H}$ given by
\begin{eqnarray*}
&&\esc{,}:H\otimes\hat{H}\rightarrow \c, \;\;\;
a\otimes \vphi \mapsto \esc{a,\vphi}:=\vphi(a),\\
&&\esc{,}:\hat{H}\otimes H\rightarrow \c, \;\;\; \vphi\otimes
a\mapsto \esc{\vphi,a}:=\vphi(a),
\end{eqnarray*}
that give a structure of \dtext{dual pairing of Hopf algebras}
between $H$ and $\hat{H}$. Associated to this pairing we have the
natural actions
\begin{eqnarray*}
&&\vartriangleright: H\otimes \hat{H}\rightarrow \hat{H}, \;\;\;
a\otimes \vphi \mapsto \vphi_{1}\esc{a,\vphi_{2}},\\
&&\vartriangleleft: \hat{H}\otimes H\rightarrow \hat{H}, \;\;\;
\vphi \otimes a\mapsto \esc{\vphi_{1},a}\vphi_{2}.
\end{eqnarray*}

Now, for every $i\in \z$, let us take $A_i:=\hat{H}$ if $i$ is odd
and $A_i:=H$ if $i$ is even, and define the maps:

\begin{eqnarray*}
    R_{2k\,2k+1}: A_{2k+1}\otimes A_{2k} & \lto  & A_{2k}\otimes A_{2k+1},\\
    \vphi\otimes a & \lmto & (\vphi_{1}\vartriangleright a)\otimes \vphi_{2} =
    a_{1}\esc{a_{2},\vphi_{1}} \otimes  \vphi_{2},\\
    R_{2k-1\,2k}: A_{2k}\otimes A_{2k-1} & \lto & A_{2k-1}\otimes A_{2k},\\
    a \otimes \vphi& \lmto & (a_{1}\vartriangleright \vphi) \otimes  a_{2} =
    \vphi_{1}\esc{\vphi_{2},a_{1}}\otimes  a_{2},\\
    R_{ij}:A_j\otimes A_i & \lto & A_i\otimes A_j,\\
    a\otimes b & \lmto & b\otimes a,\quad \text{whenever
    $j-i>2$.}
\end{eqnarray*}
As all the maps $R_{ij}$ are either usual flips or the maps induced
by a module algebra action, it is clear that all of them are
twisting maps. Furthermore, it is easy to check that they satisfy
condition \eqref{starcondition}, so they define an involution on
every twisted tensor product. Let us now check that these maps are
compatible. More precisely, let $i<j<k$, and consider the three maps
$R_{ij}$, $R_{jk}$, and $R_{ik}$, and let us check that they satisfy
the hexagon equation. We have to distinguish among several cases:
\begin{itemize}
    \item If both $\abs{j-i}, \abs{k-j}\geq 2$, all three maps
    are just usual flips, and thus the hexagon condition is
    trivially satisfied.
    \item If $\abs{j-i}=1$, $\abs{k-j}\geq 2$, then we have that
    both $R_{ik}$ and $R_{jk}$ are usual flips, and so
    the compatibility of $R_{ij}$ with them follows from Lemma
    \ref{twoflipscondition}. The same thing happens if $\abs{k-j}=1$,
    $\abs{j-i}\geq 2$.
    \item If $j=i+1, k=i+2$, then only the map $R_{i\,i+2}$ is a
    flip. Then we face two possible situations.

    If $i=2n-1$ is odd, then, describing explicitly the maps, we have
    that
    \begin{eqnarray*}
        R_{2n-1\,2n}(a\otimes
        \vphi)=\esc{\vphi_{2},a_{1}}\vphi_{1}\otimes
        a_{2},\\
        R_{2n\,2n+1}(\vphi\otimes b) = \esc{b_{2},\vphi_{1}}b_{1}\otimes
        \vphi_{2}.
    \end{eqnarray*}
    Hence, applying the left-hand side of the hexagon equation to a
    generator $a\otimes \vphi\otimes b$ of $A_{2n+1}\otimes
    A_{2n}\otimes A_{2n-1}= H\otimes \hat{H}\otimes H$, we have
    \begin{gather*}
        (A_{2n-1}\otimes R_{2n\,2n-1})(\tau\otimes
        A_{2n})(A_{2n-1}\otimes R_{2n-1\,2n})(a\otimes b\otimes c)
         \\ = (A_{2n-1}\otimes R_{2n\,2n-1})(\tau\otimes
        A_{2n})(\esc{b_{2},\vphi_{1}}a\otimes b_{1}\otimes
        \vphi_{2})
        \\=(A_{2n-1}\otimes R_{2n\,2n-1})(\esc{b_{2},\vphi_{1}}a\otimes
        \vphi_{2}\otimes b_{1})  \\
        = \esc{b_{2},\vphi_{1}}\esc{\vphi_{3},a_{1}}b_{1}\otimes
        \vphi_{2}\otimes a_{2}.
    \end{gather*}
    On the other hand, for the right hand side we get
    \begin{gather*}
        (R_{2n-1\,2n}\otimes A_{2n+1})(A_{2n}\otimes
        \tau)(R_{2n\,2n+1}\otimes A_{2n-1})(a\otimes \vphi\otimes
        b)\\
        = (R_{2n-1\,2n}\otimes A_{2n+1})(A_{2n}\otimes
        \tau)(\esc{\vphi_{2},a_{1}}\vphi_{1}\otimes a_{1} \otimes
        b)  \\
        = (R_{2n-1\,2n}\otimes A_{2n+1})(\esc{\vphi_{2},a_{1}}\vphi_{1}
         \otimes b\otimes a_{1})  \\
        =\esc{b_{2},\vphi_{1}}\esc{\vphi_{3},a_{1}}b_{1}\otimes
        \vphi_{2}\otimes a_{2},
    \end{gather*}
    where for both expressions we are using the coassociativity of
    $\hat{H}$. This proves the hexagon condition for $i$ odd. For
    $i$ even, the proof is very similar.
\end{itemize}

Now, once proven that any three twisting maps chosen from the above
ones are compatible, we can apply the Coherence Theorem and build
any iterated twisted tensor product of these algebras. In
particular, for any $n\leq m\in \z$ we may define the algebras
    \[
            A_{n, m}:=A_n\otimes_{R_{n\,n+1}} A_{n+1}
         \otimes \dotsb \otimes_{R_{m-1\, m}} A_m.
    \]
It is easy to see that if $n'\leq n$ and $m\leq m'$, then
$A_{n,m}\subseteq A_{n',m'}$ and hence the inclusions give us a
direct system of algebras $\{A_{n,m}\}_{n,m\in\z}$, being its direct
limit $\dirlim A_{n,m}$ precisely the observable algebra
$\mathcal{A}$ defined in \cite{Nill97a}. Furthermore, as the action
that defines the twisting map is a $\ast$--Hopf algebra action, we
have an involution defined on any of these products, and all the
involved algebras being of finite dimension, we have no problem
involving nuclearity nor completeness, and henceforth all the
algebras $A_{n,m}$ are well defined, finite dimensional
$C^\ast$--algebras (a fact that was proven in \cite{Nill97a} using
representations of these algebras on some Hilbert spaces). In
particular, we get a new proof of the fact that the algebra
$\mathcal{A}$ is an AF--algebra.


\section{Invariance under twisting}\label{sec:invar}
\setcounter{equation}{0}
We begin this section with a result which does not involve a twisted
tensor product of algebras and which is of independent interest.
\begin{proposition} \label{prop1}
Let $A, B$ be two algebras and $R:B\ot A\rightarrow A\ot B$ a linear
map, with notation $R(b\ot a)=a_R\ot b_R$, for all $a\in A$ and
$b\in B$. Assume that we are given two linear maps, $\mu :B\ot
A\rightarrow A$, $\mu (b\ot a)=b\cdot a$, and $\rho :A\rightarrow
A\ot B$, $\rho (a)= a_{(0)}\ot a_{(1)}$, and denote
$a*a':=a_{(0)}(a_{(1)}\cdot a')$, for all $a, a'\in A$. Assume that
the following conditions are satisfied: $\rho (1)=1\ot 1$, $1\cdot
a=a$, $a_{(0)}(a_{(1)}\cdot 1)=a$, for all $a\in A$, and
\begin{eqnarray}
&&b\cdot (a*a')=a_{(0)_R}(b_Ra_{(1)}\cdot a'), \label{a1}\\
&&\rho (a*a')=a_{(0)}a'_{(0)_R}\ot a_{(1)_R}a'_{(1)}, \label{a2}
\end{eqnarray}
for all $a, a'\in A$ and $b\in B$. Then $(A, *, 1)$ is an
associative unital algebra, denoted in what follows by $A^d$.
\end{proposition}
\begin{proof}
Obviously $1$ is the unit, so we only prove the associativity of
$*$; we compute:
\begin{eqnarray*}
(a*a')*a''&=&(a*a')_{(0)}((a*a')_{(1)}\cdot a'')\\
&\overset{\eqref{a2}}{=}&a_{(0)}a'_{(0)_R}(a_{(1)_R}a'_{(1)}\cdot
a''),
\end{eqnarray*}
\begin{eqnarray*}
a*(a'*a'')&=&a_{(0)}(a_{(1)}\cdot (a'*a''))\\
&\overset{\eqref{a1}}{=}&a_{(0)}a'_{(0)_R}(a_{(1)_R}a'_{(1)}\cdot
a''),
\end{eqnarray*}
and we see that the two terms are equal. \qed\end{proof}
\begin{remark}
The datum in Proposition \ref{prop1} is a generalization of the
left-right version of a so-called left twisting datum in
\cite{Ferrer99a}, which is obtained if $B$ is a bialgebra and the
map $R$ is given by $R(b\ot a)=b_1\cdot a\ot b_2$.
\end{remark}

As a consequence of Proposition \ref{prop1} we can obtain the
following result from \cite{Beattie96a}:

\begin{corollary} (\cite{Beattie96a}) Let $H$ be a bialgebra and $A$ a
right $H$-comodule algebra with comodule structure $A\rightarrow
A\ot H$, $a\mapsto a_{(0)}\ot a_{(1)}$, together with a linear map
$H\ot A\rightarrow A$, $h\ot a \mapsto h\cdot a$, satisfying $1\cdot
a=a$, $h\cdot 1=\varepsilon (h)1$, for all $h\in H$, $a\in A$, and
\begin{eqnarray}
&&(h_2\cdot a)_{(0)}\ot h_1(h_2\cdot a)_{(1)}=h_1\cdot a_{(0)}\ot
h_2a_{(1)}, \label{bea1} \\
&&h\cdot (a*a')=(h_1\cdot a_{(0)})(h_2a_{(1)}\cdot a'), \label{bea2}
\end{eqnarray}
where we denoted $a*a'=a_{(0)}(a_{(1)}\cdot a')$. Then $(A, *, 1)$
is an associative algebra.
\end{corollary}

\begin{proof}
We take $B=H$ and $R:H\ot A\rightarrow A\ot H$, $R(h\ot a)=h_1\cdot
a\ot h_2$. Then (\ref{a1}) is exactly (\ref{bea2}) and (\ref{a2}) is
an easy consequence of (\ref{bea1}) and of the fact that $A$ is a
comodule algebra. \qed\end{proof}
\begin{theorem} \label{main}
Assume that the hypotheses of Proposition \ref{prop1} are satisfied,
such that moreover $R$ is a twisting map. Assume also that we are
given a linear map $\lambda :A\rightarrow A\ot B$, with notation
$\lambda (a)=a_{[0]}\ot a_{[1]}$, such that $\lambda (1)=1\ot 1$ and
the following relations hold:
\begin{eqnarray}
&&\lambda (aa')=a_{[0]}*(a'_R)_{[0]}\ot (a'_R)_{[1]}(a_{[1]})_R,
\label{a3} \\
&&a_{(0)_{[0]}}\ot a_{(0)_{[1]}}a_{(1)}=a\ot 1, \label{a4}\\
&&a_{[0]_{(0)}}\ot a_{[0]_{(1)}}a_{[1]}=a\ot 1, \label{a5}
\end{eqnarray}
for all $a, a'\in A$. Define the map
\begin{eqnarray}
&&R^d:B\ot A^d\rightarrow A^d\ot B, \;\;R^d(b\ot
a)=(a_{(0)_R})_{[0]}\ot (a_{(0)_R})_{[1]}b_Ra_{(1)}. \label{stanga}
\end{eqnarray}
Then $R^d$ is a twisting map and we have an algebra isomorphism
\begin{eqnarray*}
&&A^d\ot _{R^d}B\simeq A\ot _RB, \;\;a\ot b\mapsto a_{(0)}\ot
a_{(1)}b.
\end{eqnarray*}
\end{theorem}

\begin{proof}
It is easy to see that $R^d$ satisfies \eqref{unitconditions}. We
check (\ref{tw4}) for $R^d$; we compute (denoting also
$R=r=\mathcal{R}=\overline{r}$ copies of $R$):
\begin{eqnarray*}
(a*a')_{R^d}\ot b_{R^d}&=&((a*a')_{(0)_R})_{[0]}\ot
((a*a')_{(0)_R})_{[1]}
b_R(a*a')_{(1)}\\
&\overset{\eqref{a2}}{=}&((a_{(0)}a'_{(0)_r})_R)_{[0]}\ot
((a_{(0)}a'_{(0)_r})_R)_{[1]}b_Ra_{(1)_r}a'_{(1)}\\
&\overset{\eqref{tw4}}{=}&(a_{(0)_R}(a'_{(0)_r})_\mathcal{R})_{[0]}\ot
(a_{(0)_R}(a'_{(0)_r})_\mathcal{R})_{[1]}(b_R)_\mathcal{R}a_{(1)_r}a'_{(1)}\\
&\overset{\eqref{a3}}{=}&(a_{(0)_R})_{[0]}*(((a'_{(0)_r})_\mathcal{R})_{\overline{r}}
)_{[0]}\ot (((a'_{(0)_r})_\mathcal{R})_{\overline{r}})_{[1]}
((a_{(0)_R})_{[1]})_{\overline{r}}\\
&&(b_R)_\mathcal{R}a_{(1)_r}a'_{(1)},
\end{eqnarray*}
\begin{eqnarray*}
a_{R^d}*a'_{r^d}\ot (b_{R^d})_{r^d}&=&(a_{(0)_R})_{[0]}*a'_{r^d}\ot
((a_{(0)_R})_{[1]}b_Ra_{(1)})_{r^d}\\
&=&(a_{(0)_R})_{[0]}*(a'_{(0)_r})_{[0]}\ot (a'_{(0)_r})_{[1]}
((a_{(0)_R})_{[1]}b_Ra_{(1)})_ra'_{(1)}\\
&\overset{\eqref{tw5}}{=}&(a_{(0)_R})_{[0]}*(((a'_{(0)_r})_\mathcal{R})_{\overline{r}}
)_{[0]}\ot (((a'_{(0)_r})_\mathcal{R})_{\overline{r}})_{[1]}
((a_{(0)_R})_{[1]})_{\overline{r}}\\
&&(b_R)_\mathcal{R}a_{(1)_r}a'_{(1)},
\end{eqnarray*}
and we see that the two terms coincide. Now we check (\ref{tw5}) for
$R^d$; we compute:
\begin{eqnarray*}
a_{R^d}\ot (bb')_{R^d}&=&(a_{(0)_R})_{[0]}\ot
(a_{(0)_R})_{[1]}(bb')_R
a_{(1)}\\
&\overset{\eqref{tw5}}{=}&((a_{(0)_R})_r)_{[0]}\ot
((a_{(0)_R})_r)_{[1]} b_rb'_Ra_{(1)},
\end{eqnarray*}
\begin{eqnarray*}
(a_{R^d})_{r^d}\ot b_{r^d}b'_{R^d}&=&((a_{(0)_R})_{[0]})_{r^d}\ot
b_{r^d}(a_{(0)_R})_{[1]}b'_Ra_{(1)}\\
&=&((((a_{(0)_R})_{[0]})_{(0)})_r)_{[0]}\ot
((((a_{(0)_R})_{[0]})_{(0)})_r)_{[1]}b_r\\
&&((a_{(0)_R})_{[0]})_{(1)}(a_{(0)_R})_{[1]}b'_Ra_{(1)}\\
&\overset{\eqref{a5}}{=}&((a_{(0)_R})_r)_{[0]}\ot
((a_{(0)_R})_r)_{[1]} b_rb'_Ra_{(1)},
\end{eqnarray*}
and we see that the two terms coincide, hence indeed $R^d$ is a
twisting map. We prove now that the map $\varphi :A^d\ot
_{R^d}B\rightarrow A\ot _RB$, $\varphi (a\ot b)=a_{(0)}\ot
a_{(1)}b$, is an algebra isomorphism. First, using (\ref{a4}) and
(\ref{a5}), it is easy to see that $\varphi $ is bijective, with
inverse given by $a\ot b\mapsto a_{[0]}\ot a_{[1]}b$. It is obvious
that $\varphi (1\ot 1)=1\ot 1$, so we only have to prove that
$\varphi $ is multiplicative. We compute:
\begin{eqnarray*}
\varphi ((a\ot b)(a'\ot b'))&=&\varphi (a*a'_{R^d}\ot b_{R^d}b')\\
&=&\varphi (a*(a'_{(0)_R})_{[0]}\ot (a'_{(0)_R})_{[1]}b_Ra'_{(1)}b')\\
&=&(a*(a'_{(0)_R})_{[0]})_{(0)}\ot (a*(a'_{(0)_R})_{[0]})_{(1)}
(a'_{(0)_R})_{[1]}b_Ra'_{(1)}b'\\
&\overset{\eqref{a2}}{=}&a_{(0)}(((a'_{(0)_R})_{[0]})_{(0)})_r\ot
a_{(1)_r}((a'_{(0)_R})_{[0]})_{(1)}(a'_{(0)_R})_{[1]}b_Ra'_{(1)}b'\\
&\overset{\eqref{a5}}{=}&a_{(0)}(a'_{(0)_R})_r\ot
a_{(1)_r}b_Ra'_{(1)}b',
\end{eqnarray*}
\begin{eqnarray*}
\varphi (a\ot b)\varphi (a'\ot b')&=&(a_{(0)}\ot a_{(1)}b)
(a'_{(0)}\ot a'_{(1)}b')\\
&=&a_{(0)}a'_{(0)_R}\ot (a_{(1)}b)_Ra'_{(1)}b'\\
&\overset{\eqref{tw5}}{=}&a_{(0)}(a'_{(0)_R})_r\ot
a_{(1)_r}b_Ra'_{(1)}b',
\end{eqnarray*}
finishing the proof. \qed\end{proof}

Let $H$ be a bialgebra and $F\in H\ot H$ a 2-cocycle, that is $F$ is
invertible and satisfies
\begin{eqnarray*}
&&(\varepsilon \ot id)(F)=(id \ot \varepsilon )(F)=1, \\
&&(1\ot F)(id\ot \Delta )(F)=(F\ot 1)(\Delta \ot id)(F).
\end{eqnarray*}
We denote $F=F^1\ot F^2$ and $F^{-1}=G^1\ot G^2$. We denote by $H_F$
the Drinfeld twist of $H$, which is a bialgebra having the same
algebra structure as $H$ and comultiplication given by $\Delta
_F(h)= F\Delta (h)F^{-1}$, for all $h\in H$.

If $A$ is a left $H$-module algebra (with $H$-action denoted by
$h\ot a\mapsto h\cdot a$), the invariance under twisting of the
smash product $A\# H$ is the following result (see \cite{Majid97a},
\cite{Bulacu00a}). Define a new multiplication on $A$, by
$a*a'=(G^1\cdot a)(G^2\cdot a')$, for all $a, a'\in A$, and denote
by $A_{F^{-1}}$ the new structure; then $A_{F^{-1}}$ is a left
$H_F$-module algebra (with the same action as for $A$) and we have
an algebra isomorphism
\begin{eqnarray}
&&A_{F^{-1}}\# H_F\simeq A\# H, \;\;a\# h\mapsto G^1\cdot a\# G^2h.
\label{izonoi}
\end{eqnarray}

We prove that this result is a particular case of Theorem
\ref{main}.

We take $B=H$ and $R:H\ot A\rightarrow A\ot H$, $R(h\ot a)=h_1\cdot
a\ot h_2$, hence $A\ot _RB=A\# H$. Define the maps
\begin{eqnarray*}
&&\mu :H\ot A\rightarrow A, \;\;\mu (h\ot a)=h\cdot a, \\
&&\rho :A\rightarrow A\ot H, \;\;\rho (a)=G^1\cdot a\ot G^2,\\
&&\lambda :A\rightarrow A\ot H, \;\;\lambda (a)=F^1\cdot a\ot F^2,
\end{eqnarray*}
hence the corresponding product $*$ on $A$ is given by
\begin{eqnarray*}
&&a*a'=a_{(0)}(a_{(1)}\cdot a')=(G^1\cdot a)(G^2\cdot a'),
\end{eqnarray*}
which is exactly the product of $A_{F^{-1}}$. One can check, by
direct computation, that all the necessary conditions for applying
Theorem \ref{main} are satisfied, hence we have the twisting map
$R^d:H\ot A_{F^{-1}}\rightarrow A_{F^{-1}}\ot H$, which looks as
follows:
\begin{eqnarray*}
R^d(h\ot a)&=&(a_{(0)_R})_{[0]}\ot (a_{(0)_R})_{[1]}h_Ra_{(1)}\\
&=&(h_1\cdot a_{(0)})_{[0]}\ot (h_1\cdot a_{(0)})_{[1]}h_2a_{(1)}\\
&=&(h_1G^1\cdot a)_{[0]}\ot (h_1G^1\cdot a)_{[1]}h_2G^2\\
&=&F^1h_1G^1\cdot a\ot F^2h_2G^2\\
&=&h_{(1)}\cdot a\ot h_{(2)},
\end{eqnarray*}
where we denoted by $\Delta _F(h)=h_{(1)}\ot h_{(2)}$ the
comultiplication of $H_F$. Hence, we obtain that $A^d\ot
_{R^d}B=A_{F^{-1}}\ot _{R^d}H= A_{F^{-1}}\# H_F$, and it is obvious
that the isomorphism $A^d\ot _{R^d}B\simeq A\ot _RB$ provided by
Theorem \ref{main} coincides with the one given by (\ref{izonoi}).

Let $H$ be a finite dimensional Hopf algebra with antipode $S$. As
before, we work with the realization of the Drinfeld double on $H^{*
cop}\ot H$. A well-known theorem of Majid (see \cite{Majid91a})
asserts that if $(H, r)$ is quasitriangular then the Drinfeld double
of $H$ is isomorphic to an ordinary smash product. More explicitly,
for the realization of $D(H)$ we work with, the isomorphism is given
as follows. First, we have a left $H$-module algebra structure on
$H^*$, denoted by $\underline{H}^*$, given by (we denote $r=r^1\ot
r^2$):
\begin{eqnarray*}
&&h\cdot \varphi =h_1\rightharpoonup \varphi \leftharpoonup S^{-1}(h_2), \\
&&\varphi *\varphi '=(\varphi \leftharpoonup S^{-1}(r^1))(r^2_1
\rightharpoonup \varphi '\leftharpoonup S^{-1}(r^2_2)),
\end{eqnarray*}
for all $h\in H$ and $\varphi , \varphi '\in H^*$, and then we have
an algebra isomorphism
\begin{eqnarray}
&&\underline{H}^*\# H\simeq D(H),\;\;\;\varphi \# h\mapsto \varphi
\leftharpoonup S^{-1}(r^1)\ot r^2h. \label{izomaj}
\end{eqnarray}

We prove now that this result is also a particular case of Theorem
\ref{main}.

We take $A=H^*$, with its ordinary algebra structure, $B=H$, and
$R:H\ot H^*\rightarrow H^*\ot H$, $R(h\ot \varphi
)=h_1\rightharpoonup \varphi \leftharpoonup S^{-1}(h_3)\ot h_2$,
hence $A\ot _RB=D(H)$. Denoting $r^{-1}=u^1\ot u^2$, we define the
maps:
\begin{eqnarray*}
&&\mu :H\ot H^*\rightarrow H^*, \;\;\;\mu (h\ot \varphi )=h\cdot
\varphi
=h_1\rightharpoonup \varphi \leftharpoonup S^{-1}(h_2), \\
&&\rho :H^*\rightarrow H^*\ot H, \;\;\;\rho (\varphi )=\varphi
\leftharpoonup S^{-1}(r^1)\ot r^2, \\
&&\lambda :H^*\rightarrow H^*\ot H, \;\;\;\lambda (\varphi )=\varphi
\leftharpoonup S^{-1}(u^1)\ot u^2,
\end{eqnarray*}
hence the corresponding product $*$ on $H^*$ is given by
\begin{eqnarray*}
\varphi *\varphi '&=&\varphi _{(0)}(\varphi _{(1)}\cdot \varphi ')\\
&=&(\varphi \leftharpoonup S^{-1}(r^1))(r^2\cdot \varphi ')\\
&=&(\varphi \leftharpoonup S^{-1}(r^1))(r^2_1\rightharpoonup \varphi
' \leftharpoonup S^{-1}(r^2_2)),
\end{eqnarray*}
which is exactly the product of $\underline{H}^*$. One can check, by
direct computation, that all the necessary conditions for applying
Theorem \ref{main} are satisfied, hence we have the twisting map
$R^d:H\ot \underline{H}^*\rightarrow \underline{H}^*\ot H$, which
looks as follows:
\begin{eqnarray*}
R^d(h\ot \varphi )&=&(\varphi _{(0)_R})_{[0]}\ot (\varphi
_{(0)_R})_{[1]}
h_R\varphi _{(1)}\\
&=&\varphi _{(0)_R}\leftharpoonup S^{-1}(u^1)\ot u^2h_R\varphi _{(1)}\\
&=&(\varphi \leftharpoonup S^{-1}(r^1))_R\leftharpoonup
S^{-1}(u^1)\ot
u^2h_Rr^2\\
&=&h_1\rightharpoonup \varphi \leftharpoonup S^{-1}(r^1)S^{-1}(h_3)
S^{-1}(u^1)\ot u^2h_2r^2\\
&=&h_1\rightharpoonup \varphi \leftharpoonup S^{-1}(u^1h_3r^1)\ot
u^2h_2r^2\\
&=&h_1\rightharpoonup \varphi \leftharpoonup S^{-1}(h_2)\ot h_3\\
&=&h_1\cdot \varphi \ot h_2
\end{eqnarray*}
(for the sixth equality we used the fact that $\Delta ^{cop}(h)r=
r\Delta (h)$), hence we obtain that $A^d\ot _{R^d}B=
\underline{H}^*\ot _{R^d}H=\underline{H}^*\# H$, and it is obvious
that the isomorphism $A^d\ot _{R^d}B\simeq A\ot _RB$ provided by
Theorem \ref{main} coincides with the one given by (\ref{izomaj}).

Proposition \ref{prop1} and Theorem \ref{main} admit right-left
versions, whose proofs are similar to the left-right versions above
and therefore will be omitted:
\begin{proposition} \label{prop2}
Let $B, C$ be two algebras and $R:C\ot B\rightarrow B\ot C$ a linear
map, with notation $R(c\ot b)=b_R\ot c_R$, for all $b\in B$ and
$c\in C$. Assume that we are given two linear maps, $\nu :C\ot
B\rightarrow C$, $\nu (c\ot b)=c\cdot b$, and $\theta :C\rightarrow
B\ot C$, $\theta (c)= c_{<-1>}\ot c_{<0>}$, and denote $c*c'=(c\cdot
c'_{<-1>})c'_{<0>}$, for all $c, c'\in C$. Assume that the following
conditions are satisfied: $\theta (1)=1\ot 1$, $c\cdot 1=c$,
$(1\cdot c_{<-1>})c_{<0>}=c$, for all $c\in C$, and
\begin{eqnarray}
&&(c*c')\cdot b=(c\cdot c'_{<-1>}b_R)c'_{<0>_R}, \label{b1}\\
&&\theta (c*c')=c_{<-1>}c'_{<-1>_R}\ot c_{<0>_R}c'_{<0>}, \label{b2}
\end{eqnarray}
for all $c, c'\in C$ and $b\in B$. Then $(C, *, 1)$ is an
associative unital algebra, denoted in what follows by $^dC$.
\end{proposition}
\begin{theorem} \label{maindr}
Assume that the hypotheses of Proposition \ref{prop2} are satisfied,
such that moreover $R$ is a twisting map. Assume also that we are
given a linear map $\gamma :C\rightarrow B\ot C$, with notation
$\gamma (c)=c_{\{-1\}}\ot c_{\{0\}}$, such that $\gamma (1)=1\ot 1$
and the following relations hold:
\begin{eqnarray}
&&\gamma (cc')=c'_{\{-1\}_R}(c_R)_{\{-1\}}\ot
(c_R)_{\{0\}}*c'_{\{0\}},
\label{b3} \\
&&c_{<-1>}c_{<0>_{\{-1\}}}\ot c_{<0>_{\{0\}}}=1\ot c, \label{b4}\\
&&c_{\{-1\}}c_{\{0\}_{<-1>}}\ot c_{\{0\}_{<0>}}=1\ot c, \label{b5}
\end{eqnarray}
for all $c, c'\in C$. Define the map
\begin{eqnarray}
&&^dR:\;^dC\ot B\rightarrow B\ot \;^dC, \;\;^dR(c\ot b)=c_{<-1>}b_R
(c_{<0>_R})_{\{-1\}}\ot (c_{<0>_R})_{\{0\}}. \label{dreapta}
\end{eqnarray}
Then $^dR$ is a twisting map and we have an algebra isomorphism
\begin{eqnarray*}
&&B\ot _{^dR}\;^dC\simeq B\ot _RC, \;\;b\ot c\mapsto bc_{<-1>}\ot
c_{<0>}.
\end{eqnarray*}
\end{theorem}

A particular case of Theorem \ref{maindr} is the invariance under
twisting of the right smash product from \cite{BulacuUNa}. Namely,
let $H$ be a bialgebra, $C$ a right $H$-module algebra (with action
denoted by $c\ot h\mapsto c\cdot h$) and $F\in H\ot H$ a 2-cocycle.
The right smash product $H\# C$ has multiplication
\begin{eqnarray*}
&&(h\# c)(h'\# c')=hh'_1\# (c\cdot h'_2)c'.
\end{eqnarray*}
If we define a new multiplication on $C$, by $c*c'=(c\cdot F^1)
(c'\cdot F^2)$ and denote the new structure by $_FC$, then $_FC$
becomes a right $H_F$-module algebra and we have an algebra
isomorphism $H_F\# \;_FC\simeq H\# C$, $h\# c\mapsto hF^1\# c\cdot
F^2$, see \cite{BulacuUNa}. This result may be reobtained as a
consequence of Theorem \ref{maindr}, by taking $B=H$, $R(c\ot
h)=h_1\ot c\cdot h_2$, $\nu (c\ot h)=c\cdot h$, $\theta (c)=F^1\ot
c\cdot F^2$, $\gamma (c)=G^1\ot c\cdot G^2$, where we denoted as
before $F^{-1}=G^1\ot G^2$.

A careful look at the proof of Theorem \ref{main} shows that
actually it admits a more general form, which we record here for
further use (the same holds for Theorem \ref{maindr}).
\begin{theorem} \label{degen1}
Let $A\ot _RB$ be a twisted tensor product of algebras, and denote
the multiplication of $A$ by $a\ot a'\mapsto aa'$. Assume that on
the vector space $A$ we have one more algebra structure, denoted by
$A'$, with the same unit as $A$ and multiplication denoted by $a\ot
a'\mapsto a*a'$ (for instance, $A'$ may be $A$ itself or $A^d$ as in
Proposition \ref{prop1}). Assume that we are given two linear maps
$\rho ,\lambda :A\rightarrow A\ot B$, with notation $\rho
(a)=a_{(0)}\ot a_{(1)}$ and $\lambda (a)=a_{[0]}\ot a_{[1]}$, such
that $\rho $ is an algebra map from $A'$ to $A\ot _RB$, $\lambda
(1)=1\ot 1$ and relations (\ref{a3}), (\ref{a4}), (\ref{a5}) are
satisfied. Then the map
\begin{eqnarray}
&&R':B\ot A'\rightarrow A'\ot B, \;\;R'(b\ot a)=(a_{(0)_R})_{[0]}\ot
(a_{(0)_R})_{[1]}b_Ra_{(1)}, \label{degstanga}
\end{eqnarray}
is a twisting map and we have an algebra isomorphism
\begin{eqnarray*}
&&A'\ot _{R'}B\simeq A\ot _RB, \;\;a\ot b\mapsto a_{(0)}\ot
a_{(1)}b.
\end{eqnarray*}
\end{theorem}
\begin{theorem} \label{degen2}
Let $B\ot _RC$ be a twisted tensor product of algebras and denote
the multiplication of $C$ by $c\ot c'\mapsto cc'$. Assume that on
the vector space $C$ we have one more algebra structure, denoted by
$C'$, with the same unit as $C$ and multiplication denoted by $c\ot
c'\mapsto c*c'$ (for instance, $C'$ may be $C$ itself or $^dC$ as in
Proposition \ref{prop2}). Assume that we are given two linear maps
$\theta ,\gamma :C\rightarrow B\ot C$, with notation $\theta
(c)=c_{<-1>}\ot c_{<0>}$ and $\gamma (c)=c_{\{-1\}}\ot c_{\{0\}}$,
such that $\theta $ is an algebra map from $C'$ to $B\ot _RC$,
$\gamma (1)=1\ot 1$ and relations (\ref{b3}), (\ref{b4}), (\ref{b5})
are satisfied. Then the map
\begin{eqnarray}
&&R':\;C'\ot B\rightarrow B\ot C', \;\;R'(c\ot b)=c_{<-1>}b_R
(c_{<0>_R})_{\{-1\}}\ot (c_{<0>_R})_{\{0\}}, \label{degdreapta}
\end{eqnarray}
is a twisting map and we have an algebra isomorphism
\begin{eqnarray*}
&&B\ot _{R'}C'\simeq B\ot _RC, \;\;b\ot c\mapsto bc_{<-1>}\ot
c_{<0>}.
\end{eqnarray*}
\end{theorem}


We recall the following result of G. Fiore from \cite{fiore}, in a
slightly modified (but equivalent) form. Let $H$ be a Hopf algebra with
antipode $S$ and $A$ a left $H$-module algebra. Assume that there exists an
algebra map $\varphi :A\# H\rightarrow A$ such that $\varphi (a\# 1)=a$
for all $a\in A$. Define the map
\begin{eqnarray*}
&&\theta :H\rightarrow A\otimes H, \;\;\;\theta (h)=\varphi (1\# S(h_1))
\otimes h_2.
\end{eqnarray*}
Then $\theta $ is an algebra map from $H$ to $A\# H$ and the smash
product $A\# H$ is isomorphic to the ordinary tensor product $A\otimes H$.

We prove that this result is a particular case of Theorem \ref{degen2},
with $B=A$ and $C=C'=H$ (in the notation of Theorem \ref{degen2}).

Define the map $\gamma :H\rightarrow A\otimes H$, $\gamma (h)=
\varphi (1\# h_1)\otimes h_2$, and denote as above $\theta (h)=h_{<-1>}
\otimes h_{<0>}$ and $\gamma (h)=h_{\{-1\}}\otimes h_{\{0\}}$. The
relations (\ref{b4}) and (\ref{b5}) are easy to check, so we only have to
prove (\ref{b3}) (here, the map $R:H\otimes A\rightarrow A\otimes H$ is
given by $R(h\otimes a)=h_1\cdot a \otimes h_2$). We will need the
following relation from \cite{fiore}:
\begin{eqnarray}
&&\varphi (1\# h)a=(h_1\cdot a)\varphi (1\# h_2), \label{fio}
\end{eqnarray}
for all $h\in H$, $a\in A$. Now we compute:
\begin{eqnarray*}
(h'_{\{-1\}})_R(h_R)_{\{-1\}}\otimes (h_R)_{\{0\}}h'_{\{0\}}&=&
\varphi (1\# h'_1)_R\varphi (1\# (h_R)_1)\otimes (h_R)_2h'_2\\
&=&(h_1\cdot \varphi (1\# h'_1))\varphi (1\# h_2)\otimes h_3h'_2\\
&\overset{{\rm (\ref{fio})}}{=}&\varphi (1\# h_1)\varphi (1\# h'_1)\otimes h_2h'_2\\
&=&\varphi (1\# h_1h'_1)\otimes h_2h'_2\\
&=&\gamma (hh'),
\end{eqnarray*}
hence (\ref{b3}) holds. Theorem \ref{degen2} may thus be applied, and we
get the twisting map $R'$, which looks as follows:
\begin{eqnarray*}
R'(h\otimes a)&=&h_{<-1>}a_R(h_{<0>_R})_{\{-1\}}\otimes
(h_{<0>_R})_{\{0\}}\\
&=&\varphi (1\# S(h_1))a_R(h_{2_R})_{\{-1\}}\otimes (h_{2_R})_{\{0\}}\\
&=&\varphi (1\# S(h_1))(h_2\cdot a)(h_3)_{\{-1\}}\otimes (h_3)_{\{0\}}\\
&=&\varphi (1\# S(h_1))(h_2\cdot a)\varphi (1\# h_3)\otimes h_4\\
&\overset{{\rm (\ref{fio})}}{=}&\varphi (1\# S(h_1))\varphi (1\# h_2)a\otimes h_3\\
&=&\varphi (1\# S(h_1)h_2)a\otimes h_3\\
&=&a\otimes h,
\end{eqnarray*}
so $R'$ is the usual flip, hence we obtain $A\# H\simeq A\otimes H$ as a
consequence of Theorem \ref{degen2}.
\begin{remark}
Let $H$ be a Hopf algebra, let $A$ be an algebra and $u:H\rightarrow A$
an algebra map; consider the strongly inner action of $H$ on $A$
afforded by $u$, that is $h\cdot a=u(h_1)au(S(h_2))$, for all $h\in H$,
$a\in A$. Then it is well-known (see for instance \cite{montgomery},
Example 7.3.3) that the smash product $A\# H$ is isomorphic to the
ordinary tensor product $A\otimes H$. This result is actually a particular
case of Fiore's theorem presented above (hence of Theorem \ref{degen2} too),
because one can easily see that the map $\varphi :A\# H\rightarrow A$,
$\varphi (a\# h)=au(h)$ is an algebra map satisfying $\varphi (a\# 1)=a$
for all $a\in A$.
\end{remark}


We recall now the following result from \cite{Fiore03a}, with a
different notation and in a slightly modified (but equivalent) form,
adapted to our purpose. Let $(H, r)$ be a quasitriangular Hopf
algebra, $H^+$ and $H^-$ two Hopf subalgebras of $H$ such that $r\in
H^+\ot H^-$ (we will denote $r=r^1\ot r^2=\mathcal{R}^1\ot
\mathcal{R}^2\in H^+\ot H^-$). Let $B$ be a right $H^+$-module
algebra and $C$ a right $H^-$-module algebra (actions are denoted by
$\cdot $), and consider their braided product $B\underline{\ot }C$,
which is just the twisted tensor product $B\ot _RC$, with twisting
map given by
\begin{eqnarray*}
&&R:C\ot B\rightarrow B\ot C, \;\;\;R(c\ot b)=b\cdot r^1\ot c\cdot
r^2.
\end{eqnarray*}
Assume that there exists an algebra map $\pi :H^+\# B\rightarrow B$
(where $H^+\# B$ is the right smash product recalled before) such
that $\pi (1\# b)=b$ for all $b\in B$. Define the map
\begin{eqnarray*}
&&\theta :C\rightarrow B\ot C, \;\;\;\theta (c)=\pi (r^1\# 1)\ot
c\cdot r^2.
\end{eqnarray*}
Then $\theta $ is an algebra map from $C$ to $B\underline{\ot }C$
and the braided tensor product $B\underline{\ot }C$ is isomorphic to
the ordinary tensor product $B\ot C$ (hence the existence of $\pi $
allows to ``unbraid'' the braided tensor product; many examples
where this happens may be found in \cite{Fiore03a}, especially
coming from quantum groups).

We prove now that this result is a particular case of Theorem
\ref{degen2}, with $C'=C$ (in the notation of Theorem \ref{degen2}).

We first need to recall the axioms of a quasitriangular structure:
\begin{eqnarray}
&&(\Delta \ot id)(r)=r_{13}r_{23}, \label{qt1} \\
&&(id \ot \Delta )(r)=r_{13}r_{12}, \label{qt2} \\
&&\Delta ^{cop} (h)r=r\Delta (h), \;\;\;\forall \;h\in H.
\end{eqnarray}
Define the map $\gamma :C\rightarrow B\ot C$, $\gamma (c)=\pi (u^1\#
1)\ot c\cdot u^2$, where we denote $r^{-1}=u^1\ot u^2=U^1\ot U^2\in
H^+\ot H^-$. Denote as above $\theta (c)=c_{<-1>}\ot c_{<0>}$ and
$\gamma (c)=c_{\{-1\}}\ot c_{\{0\}}$. The relations (\ref{b4}) and
(\ref{b5}) are easy to check, hence we only have to prove (\ref{b3})
(here, we recall, $*$ coincides with the multiplication of $C$). We
first record the relation:
\begin{eqnarray}
&&c_{\{-1\}}b\ot c_{\{0\}}=b_R(c_R)_{\{-1\}}\ot (c_R)_{\{0\}},
\;\;\; \forall \;b\in B, \;c\in C, \label{ajut}
\end{eqnarray}
which can be proved as follows:
\begin{eqnarray*}
b_R(c_R)_{\{-1\}}\ot (c_R)_{\{0\}}&=&
(b\cdot r^1)(c\cdot r^2)_{\{-1\}}\ot (c\cdot r^2)_{\{0\}}\\
&=&(b\cdot r^1)\pi (u^1\# 1)\ot c\cdot r^2u^2\\
&=&\pi (1\# b\cdot r^1)\pi (u^1\# 1)\ot c\cdot r^2u^2\\
&=&\pi ((1\# b\cdot r^1)(u^1\# 1))\ot c\cdot r^2u^2\\
&=&\pi (u^1_1\# b\cdot r^1u^1_2)\ot c\cdot r^2u^2\\
&\overset{\eqref{qt1}}{=}&\pi (U^1\# b\cdot r^1u^1)\ot c\cdot r^2u^2U^2\\
&=&\pi (U^1\# b)\ot c\cdot U^2\\
&=&\pi (U^1\# 1)\pi (1\# b)\ot c\cdot U^2\\
&=&\pi (U^1\# 1)b\ot c\cdot U^2\\
&=&c_{\{-1\}}b\ot c_{\{0\}}.
\end{eqnarray*}
Now we compute:
\begin{eqnarray*}
\gamma (cc')&=&\pi (u^1\# 1)\ot (c\cdot u_1^2)(c'\cdot u_2^2)\\
&\overset{\eqref{qt2}}{=}&\pi (u^1U^1\# 1)\ot (c\cdot u^2)(c'\cdot U^2)\\
&=&\pi (u^1\# 1)\pi (U^1\# 1)\ot (c\cdot u^2)(c'\cdot U^2)\\
&=&c_{\{-1\}}c'_{\{-1\}}\ot c_{\{0\}}c'_{\{0\}}\\
&\overset{\eqref{ajut}}{=}&c'_{{\{-1\}}_R}(c_R)_{\{-1\}} \ot
(c_R)_{\{0\}}c'_{\{0\}},
\end{eqnarray*}
hence (\ref{b3}) holds. Theorem \ref{degen2} may thus be applied,
and we get the twisting map $R'$, which looks as follows:
\begin{eqnarray*}
R'(c\ot b)&=&c_{<-1>}b_R(c_{<0>_R})_{\{-1\}}\ot (c_{<0>_R})_{\{0\}}\\
&=&\pi (r^1\# 1)b_R((c\cdot r^2)_R)_{\{-1\}}\ot ((c\cdot r^2)_R)_{\{0\}}\\
&=&\pi (r^1\# 1)(b\cdot \mathcal{R}^1)(c\cdot
r^2\mathcal{R}^2)_{\{-1\}}\ot
(c\cdot r^2\mathcal{R}^2)_{\{0\}}\\
&=&\pi (r^1\# 1)\pi (1\# b\cdot \mathcal{R}^1)\pi (u^1\# 1)\ot
c\cdot r^2\mathcal{R}^2u^2\\
&=&\pi ((r^1\# 1)(1\# b\cdot \mathcal{R}^1)(u^1\# 1))\ot c\cdot
r^2\mathcal{R}^2u^2\\
&=&\pi (r^1u_1^1\# b\cdot \mathcal{R}^1u_2^1)\ot c\cdot r^2\mathcal{R}^2u^2\\
&\overset{\eqref{qt1}}{=}&\pi (r^1U^1\# b\cdot \mathcal{R}^1u^1)\ot
c\cdot r^2
\mathcal{R}^2u^2U^2\\
&=&\pi (1\# b)\ot c\\
&=&b\ot c,
\end{eqnarray*}
so $R'$ is the usual flip, hence we obtain $B\underline{\ot }C\simeq
B\ot C$ as a consequence of Theorem \ref{degen2}.

A natural question that arises is to see whether Theorems \ref{main}
and \ref{maindr} can be combined, namely, if $(A, B, C, R_1, R_2,
R_3)$ are as in Theorem \ref{itertwisting} and we have a datum as in
Theorem \ref{main} between $A$ and $B$ and a datum as in Theorem
\ref{maindr} between $B$ and $C$, under what conditions it follows
that $(A^d, B, \;^dC, R_1^d,
$$\;^dR_2, R_3)$ satisfy again the hypotheses of Theorem
\ref{itertwisting}.

Our first remark is that this does {\it not} happen in general, a
counterexample may be obtained as follows. Take $B=H$ a bialgebra,
$A$ a left $H$-module algebra, $C$ a right $H$-module algebra and
$F\in H\ot H$ a 2-cocycle. Here $R_1(h\ot a)=h_1\cdot a\ot h_2$,
$R_2(c\ot h)=h_1\ot c\cdot h_2$ and $R_3=\tau _{CA}$, the usual
flip, hence $A\ot _{R_1}H\ot _{R_2}C=A\# H\# C$, the two-sided smash
product. We consider the datum between $A$ and $H$ that allows us to
define $A_{F^{-1}}\# H_F$, hence $R_1^d(h\ot a)=F^1h_1G^1\cdot a\ot
F^2h_2G^2$, and the trivial datum between $H$ and $C$. One can see
that in general $(R_1^d, R_2, R_3)$ do {\it not} satisfy the hexagon
condition.

Hence, the best we can do is to find sufficient conditions on the
initial data ensuring that $(R_1^d, \;^dR_2, R_3)$ satisfy the
hexagon condition. This is achieved in the next result. Note that
the conditions we found are not the most general one can imagine (in
particular, we need to assume that $R_3$ is the flip), but they are
general enough to include as a particular case the invariance under
twisting of the two-sided smash product from \cite{BulacuUNa}, which
was our guiding example for this result.

\begin{theorem} \label{tare}
Let $(A, B, C, R_1, R_2, R_3)$ be as in Theorem \ref{itertwisting},
with $R_3=\tau _{CA}$, the usual flip. Assume that we have a datum
between $A$ and $B$ as in Theorem \ref{main} and a datum between $B$
and $C$ as in Theorem \ref{maindr}, with notation as in these
results. Assume also that the following compatibility conditions
hold:
\begin{eqnarray}
&&a_{(0)}\ot a_{(1)_{R_2}}(c_{R_2})_{\{-1\}}\ot (c_{R_2})_{\{0\}}=
a_{(0)_{R_1}}\ot c_{\{-1\}_{R_1}}a_{(1)}\ot c_{\{0\}}, \label{tare1} \\
&&a_{[0]}\ot c_{<-1>}a_{[1]_{R_2}}\ot c_{<0>_{R_2}}=(a_{R_1})_{[0]}
\ot (a_{R_1})_{[1]}c_{<-1>_{R_1}}\ot c_{<0>}, \label{tare2} \\
&&a_{(0)_{R_1}}\ot c_{<-1>_{R_1}}a_{(1)_{R_2}}\ot c_{<0>_{R_2}}=
a_{(0)}\ot a_{(1)}c_{<-1>}\ot c_{<0>}, \label{tare3}
\end{eqnarray}
for all $a\in A$, $c\in C$. Then $(A^d, B, \;^dC, R_1^d, $$\;^dR_2,
R_3)$ satisfy also the hypotheses of Theorem \ref{itertwisting}, and
we have an algebra isomorphism
\begin{eqnarray*}
&&A^d\ot _{R_1^d}B\ot _{\;^dR_2}\;^dC\simeq A\ot _{R_1}B\ot _{R_2}C,
\;\;\; a\ot b\ot c\mapsto a_{(0)}\ot a_{(1)}bc_{<-1>}\ot c_{<0>}.
\end{eqnarray*}
\end{theorem}

\begin{proof}
We prove the hexagon condition for $(R_1^d, \;^dR_2, R_3)$;
we compute:\\[2mm]
${\;\;}$ $(A\ot \;^dR_2)\circ (R_3\ot B)\circ (C\ot R_1^d)(c\ot b\ot
a)$
\begin{eqnarray*}
&\overset{\eqref{stanga}, \eqref{dreapta}}{=}&
(a_{(0)_{R_1}})_{[0]}\ot c_{<-1>}((a_{(0)_{R_1}})_{[1]}b_{R_1}
a_{(1)})_{R_2}
(c_{<0>_{R_2}})_{\{-1\}}\ot (c_{<0>_{R_2}})_{\{0\}}\\
&\overset{\eqref{tw4}}{=}&(a_{(0)_{R_1}})_{[0]}\ot c_{<-1>}
((a_{(0)_{R_1}})_{[1]})_{R_2}(b_{R_1})_{r_2}(a_{(1)})_{\mathcal{R}_2}
(((c_{<0>_{R_2}})_{r_2})_{\mathcal{R}_2})_{\{-1\}}\\
&&\ot (((c_{<0>_{R_2}})_{r_2})_{\mathcal{R}_2})_{\{0\}}\\
&\overset{\eqref{tare1}}{=}&(((a_{(0)})_{\mathcal{R}_1})_{R_1})_{[0]}\ot
c_{<-1>}((((a_{(0)})_{\mathcal{R}_1})_{R_1})_{[1]})_{R_2}(b_{R_1})_{r_2}\\
&&(((c_{<0>_{R_2}})_{r_2})_{\{-1\}})_{\mathcal{R}_1}a_{(1)} \ot
((c_{<0>_{R_2}})_{r_2})_{\{0\}},
\end{eqnarray*}
${\;\;}$ $(R_1^d\ot C)\circ (B\ot R_3)\circ (^dR_2\ot A)(c\ot b\ot
a)$
\begin{eqnarray*}
&\overset{\eqref{stanga}, \eqref{dreapta}}{=}&
(a_{(0)_{R_1}})_{[0]}\ot (a_{(0)_{R_1}})_{[1]}(c_{<-1>}b_{R_2}
(c_{<0>_{R_2}})_{\{-1\}})_{R_1}a_{(1)}\ot (c_{<0>_{R_2}})_{\{0\}}\\
&\overset{\eqref{tw5}}{=}&(((a_{(0)_{R_1}})_{r_1})_{\mathcal{R}_1})_{[0]}
\ot
(((a_{(0)_{R_1}})_{r_1})_{\mathcal{R}_1})_{[1]}(c_{<-1>})_{\mathcal{R}_1}
(b_{R_2})_{r_1}\\
&&((c_{<0>_{R_2}})_{\{-1\}})_{R_1}a_{(1)}\ot (c_{<0>_{R_2}})_{\{0\}}\\
&\overset{\eqref{tare2}}{=}&((a_{(0)_{R_1}})_{r_1})_{[0]}\ot
c_{<-1>}(((a_{(0)_{R_1}})_{r_1})_{[1]})_{r_2}(b_{R_2})_{r_1}
(((c_{<0>_{r_2}})_{R_2})_{\{-1\}})_{R_1}\\
&&a_{(1)}\ot ((c_{<0>_{r_2}})_{R_2})_{\{0\}},
\end{eqnarray*}
and the two terms are equal because of the hexagon condition for
$(R_1, R_2, R_3)$:
\begin{eqnarray*}
&&a_{R_1}\ot (b_{R_1})_{R_2}\ot c_{R_2}=a_{R_1}\ot (b_{R_2})_{R_1}
\ot c_{R_2}.
\end{eqnarray*}
We prove now that the map
\begin{eqnarray*}
&&\psi :A^d\ot _{R_1^d}B\ot _{\;^dR_2}\;^dC\rightarrow A\ot _{R_1}B
\ot _{R_2}C, \\
&&\psi (a\ot b\ot c)=a_{(0)}\ot a_{(1)}bc_{<-1>}\ot c_{<0>},
\end{eqnarray*}
is an algebra isomorphism. First, using (\ref{a4}), (\ref{a5}),
(\ref{b4}), (\ref{b5}), it is easy to see that $\psi $ is bijective,
with inverse given by $a\ot b\ot c\mapsto a_{[0]}\ot
a_{[1]}bc_{\{-1\}}\ot c_{\{0\}}$. We prove now that $\psi $ is
multiplicative. We compute
(using (\ref{trei})):\\[2mm]
$\psi ((a\ot b\ot c)(a'\ot b'\ot c'))$
\begin{eqnarray*}
&=&\psi (a*a'_{R_1^d}\ot b_{R_1^d}b'_{\;^dR_2}\ot c_{\;^dR_2}*c')\\
&\overset{\eqref{stanga}, \eqref{dreapta}}{=}&\psi
(a*(a'_{(0)_{R_1}})_{[0]} \ot
(a'_{(0)_{R_1}})_{[1]}b_{R_1}a'_{(1)}c_{<-1>}b'_{R_2}
(c_{<0>_{R_2}})_{\{-1\}}\\
&&\ot (c_{<0>_{R_2}})_{\{0\}}*c')\\
&=&(a*(a'_{(0)_{R_1}})_{[0]})_{(0)} \ot
(a*(a'_{(0)_{R_1}})_{[0]})_{(1)}
(a'_{(0)_{R_1}})_{[1]}b_{R_1}a'_{(1)}c_{<-1>}b'_{R_2}\\
&&(c_{<0>_{R_2}})_{\{-1\}}((c_{<0>_{R_2}})_{\{0\}}*c')_{<-1>}
\ot ((c_{<0>_{R_2}})_{\{0\}}*c')_{<0>}\\
&\overset{\eqref{a2},\eqref{b2}}{=}&a_{(0)}(((a'_{(0)_{R_1}})_{[0]})_{(0)})_{r_1}
\ot
a_{(1)_{r_1}}((a'_{(0)_{R_1}})_{[0]})_{(1)}(a'_{(0)_{R_1}})_{[1]}
b_{R_1}a'_{(1)}c_{<-1>}b'_{R_2}\\
&&(c_{<0>_{R_2}})_{\{-1\}}((c_{<0>_{R_2}})_{\{0\}})_{<-1>}
(c'_{<-1>})_{r_2}\ot (((c_{<0>_{R_2}})_{\{0\}})_{<0>})_{r_2}c'_{<0>}\\
&\overset{\eqref{a5},\eqref{b5}}{=}&a_{(0)}(a'_{(0)_{R_1}})_{r_1}\ot
a_{(1)_{r_1}}b_{R_1}a'_{(1)}c_{<-1>}b'_{R_2}c'_{<-1>_{r_2}}\ot
(c_{<0>_{R_2}})_{r_2}c'_{<0>},
\end{eqnarray*}
$\psi (a\ot b\ot c)\psi (a'\ot b'\ot c')$
\begin{eqnarray*}
&=&(a_{(0)}\ot a_{(1)}bc_{<-1>}\ot c_{<0>})(a'_{(0)}\ot a'_{(1)}b'
c'_{<-1>}\ot c'_{<0>})\\
&=&a_{(0)}a'_{(0)_{R_1}}\ot (a_{(1)}bc_{<-1>})_{R_1}(a'_{(1)}b'
c'_{<-1>})_{R_2}\ot c_{<0>_{R_2}}c'_{<0>}\\
&\overset{\eqref{tw4},\eqref{tw5}}{=}&a_{(0)}(((a'_{(0)})_{\mathcal{R}_1})_{R_1})
_{r_1}\ot
a_{(1)_{r_1}}b_{R_1}(c_{<-1>})_{\mathcal{R}_1}(a'_{(1)})_{\mathcal{R}_2}
b'_{R_2}c'_{<-1>_{r_2}}\\
&&\ot (((c_{<0>})_{\mathcal{R}_2})_{R_2})_{r_2}c'_{<0>}\\
&\overset{\eqref{tare3}}{=}&a_{(0)}(a'_{(0)_{R_1}})_{r_1}\ot
a_{(1)_{r_1}} b_{R_1}a'_{(1)}c_{<-1>}b'_{R_2}c'_{<-1>_{r_2}}\ot
(c_{<0>_{R_2}})_{r_2} c'_{<0>},
\end{eqnarray*}
and we see that the two terms are equal. \qed\end{proof}

Let now $H$ be a bialgebra, $A$ a left $H$-module algebra, $C$ a
right $H$-module algebra and $F\in H\ot H$ a 2-cocycle. Then, by
\cite{BulacuUNa}, we have an algebra isomorphism (notation as
before):
\begin{eqnarray*}
&&A_{F^{-1}}\# H_F\# \;_FC\simeq A\# H\# C, \;\;\;a\# h\# c\mapsto
G^1\cdot a\# G^2hF^1\# c\cdot F^2.
\end{eqnarray*}
One can easily see that this result is a particular case of Theorem
\ref{tare}; indeed, the relations (\ref{tare1}), (\ref{tare2}),
(\ref{tare3}) are easy consequences of the 2-cocycle condition for
$F$.

\nocite{Kassel95a} \nocite{Majid95a} \nocite{Tambara90a}
\nocite{Reshetikhin90a}

\end{document}